\documentclass[11pt]{article}
\usepackage[T1]{fontenc}
\usepackage{lmodern,amsmath,amssymb,amsthm,mathtools}
\usepackage[margin=1.05in]{geometry}
\usepackage[expansion=false]{microtype}
\usepackage{enumitem,xcolor}
\usepackage{graphicx,float,tikz}
\usetikzlibrary{calc}
\usepackage[colorlinks=true,linkcolor=blue!40!black,citecolor=blue!40!black,urlcolor=blue!40!black]{hyperref}
\setlist[enumerate]{leftmargin=2.2em,itemsep=3pt}

\newtheorem{theorem}{Theorem}
\numberwithin{theorem}{section}

\newtheorem{lemma}[theorem]{Lemma}

\theoremstyle{definition}

\newtheorem{question}[theorem]{Question}

\newtheorem{claimnum}{Claim}
\newcommand{\F}{\mathbb F_2}
\newcommand{\one}{\boldsymbol 1}
\newcommand{\R}{\mathbb R}
\newcommand{\cU}{\mathcal U}
\DeclareMathOperator{\supp}{supp}
\DeclareMathOperator{\ov}{ov}
\DeclareMathOperator{\Hess}{Hess}
\DeclareMathOperator{\tr}{tr}

\DeclareMathOperator{\conv}{conv}

\DeclareMathOperator{\im}{im}

\DeclareMathOperator{\Ind}{Ind}

\newcommand{\mc}[1]{\mathcal{#1}} 

\newcommand{\brm}[1]{\operatorname{#1}}

\title{Sharp topological overlap in every dimension}
\author{Sergey Norin\thanks{Department of Mathematics and Statistics, McGill University, Montr\'eal, Canada (\textsf{\href{mailto:sergey.norin@mcgill.ca}{sergey.norin\allowbreak@mcgill.ca}}).  Research supported by NSERC.}}
\date{\today}

\begin{document}
\maketitle

\begin{abstract}
For a positive integer $d$ let $c_d^{\mathrm{top}}$ (resp. $c_d^{\mathrm{aff}}$) be the maximum $c$ such that, for every continuous (resp. affine) map from the $n$-dimensional simplex $\Delta^n$ to $\R^d$, some point of $\R^d$ lies in the images of at least $(c-o(1))\binom{n+1}{d+1}$ of the $d$-dimensional faces, as $n\to\infty$. We show that
\[
 c_d^{\mathrm{top}}=c_d^{\mathrm{aff}}=\frac{(d+1)!}{(d+1)^{d+1}}.
\]
The second equality confirms a conjecture of Bukh, Matou\v{s}ek, and Nivasch, who established the upper bound. 

The proof was found by GPT-6 Astra, following the directions of the author. 
\end{abstract}

\section{Introduction}\label{sec:intro}

For $n\ge0$, write $[n]=\{0,\ldots,n\}$.\footnote{Note that this differs from the standard convention $[n]=\{1,\ldots,n\}$.} Let $\Delta^n$ denote the geometric $n$-dimensional simplex with vertex set $[n]$. For a continuous map $f:\Delta^n \to \R^d$, put
\begin{equation*}
 \ov(f)=\max_{y\in\R^d}\left|\left\{\sigma\in\binom{[n]}{d+1}:y\in f(\Delta(\sigma))\right\}\right|,
\end{equation*}
the largest number of $d$-dimensional faces whose images have a common point. 

Let
\begin{equation*}
 c_d^{\mathrm{top}}=\liminf_{n\to\infty}\ \min_f\frac{\ov(f)}{\binom{n+1}{d+1}},
\end{equation*}
where the minimum is over all continuous maps $f:\Delta^n\to\R^d$.

Let $c_d^{\mathrm{aff}}$ be defined analogously with the minimum taken over affine maps instead. Thus
\begin{equation*}
 c_d^{\mathrm{top}}\le c_d^{\mathrm{aff}}.
\end{equation*}
The constant $c_d^{\mathrm{aff}}$ can be equivalently defined in terms of point sets in $\R^d$. A set $S\subset\R^d$ is in \emph{general position} if no $d+1$ of its points lie on a common hyperplane; then the convex hulls of its $(d+1)$-element subsets are $\binom{|S|}{d+1}$ distinct $d$-simplices, the $d$-simplices \emph{spanned} by~$S$. For a finite $S$ in general position, let $\ov(S)$ be the largest number of $d$-simplices spanned by $S$ with a common point. Then
\begin{equation*}
 c_d^{\mathrm{aff}}=\liminf_{n\to\infty}\ \min_{|S|=n+1}\frac{\ov(S)}{\binom{n+1}{d+1}}.
\end{equation*}

\subsubsection*{Previous work}

The affine problem came first. Boros and F\"uredi~\cite{BF} proved $c_2^{\mathrm{aff}}\ge2/9$. (See also Bukh~\cite{Bukh} for a short proof.) B\'ar\'any~\cite{Barany} proved that $c_d^{\mathrm{aff}}\ge(d+1)^{-d}$ for every $d$. Wagner~\cite{Wagner} improved  the bound  to $c_d^{\mathrm{aff}}\ge(d^2+1)/(d+1)^{d+1}$. Bukh, Matou\v{s}ek, and Nivasch~\cite[Theorem~1.3]{BMN} constructed point sets showing that
\begin{equation*}
 c_d^{\mathrm{aff}}\le\theta_d:=\frac{(d+1)!}{(d+1)^{d+1}}.
\end{equation*}
Their construction is in general position, as follows from the determinant calculation in the proof of~\cite[Lemma~2.2]{BMN}. They conjectured that this is the true value. For $d=2$ their construction shows that the Boros--F\"uredi bound is sharp.

Gromov~\cite[\S1.1]{Gromov} proved the first nontrivial lower bound on the topological overlap constants in every dimension:
\begin{equation*}
 c_d^{\mathrm{top}}\ge\frac{2d}{(d+1)(d+1)!}.
\end{equation*}
For $d=2$ it equals $2/9$, so $c_2^{\mathrm{top}}=c_2^{\mathrm{aff}}=2/9$; for $d\ge3$ Gromov's bound exceeded the affine lower bounds then known.  Matou\v{s}ek and Wagner~\cite{MW} gave an exposition of the combinatorial part of Gromov's method, isolated the quantities it depends on (cofilling profiles), showed that improving them alone cannot reach $\theta_3=3/32$, and introduced pagodas to prove $c_3^{\mathrm{top}}>0.06332$~\cite[Proposition~15]{MW}. Kr\'al', Mach, and Sereni~\cite{KMS} improved one of the cofilling profiles and, combining this with the pagoda method, obtained a lower bound $c_3^{\mathrm{top}}>0.07480$. Karasev~\cite{Karasev} gave a short proof of Gromov's bound for affine maps, using only the degree of maps between spheres. 

\subsubsection*{Our results}

We determine the topological and affine overlap constants exactly in every dimension, showing that the Bukh-Matou\v{s}ek-Nivasch upper bound is tight.

\begin{theorem}\label{thm:main}
For every integer~$d\ge 1$,
\begin{equation*}
 c_d^{\mathrm{top}}=c_d^{\mathrm{aff}}=\frac{(d+1)!}{(d+1)^{d+1}}.
\end{equation*}
\end{theorem}

The second equality confirms the conjecture of Bukh, Matou\v{s}ek, and Nivasch~\cite{BMN}.

Throughout the paper we fix an integer $d\ge1$, and let $\theta= \theta_d = (d+1)!/(d+1)^{d+1}$. Since $c_d^{\mathrm{top}}\le c_d^{\mathrm{aff}}\le\theta$, Theorem~\ref{thm:main} follows from the lower bound $c_d^{\mathrm{top}}\ge\theta$, which is a consequence of the following finite statement.

\begin{theorem}\label{thm:finite}
Let
\begin{equation}\label{eq:C}
 N_d=\sum_{i=1}^{d}\binom{d+1}{i+1}(d+1-i)\left\lfloor\frac di\right\rfloor.
\end{equation}
If $n+1>N_d$, then every continuous map $f:\Delta^n\to\R^d$ satisfies
\begin{equation*}
 \ov(f)\ge\left(\frac{n+1-N_d}{d+1}\right)^{d+1}-N_d\binom{n}{d}.
\end{equation*}
\end{theorem}

\begin{proof}[Proof of Theorem~\ref{thm:main} assuming Theorem~\ref{thm:finite}]
Let $n+1>N_d$ and let $f:\Delta^n\to\R^d$ be continuous. 
Theorem~\ref{thm:finite} gives
\begin{equation*}
 \frac{\ov(f)}{\binom{n+1}{d+1}}\ge\frac{(d+1)!}{(n+1)^{d+1}}\left(\frac{n+1-N_d}{d+1}\right)^{d+1}-\frac{(d+1)N_d}{n+1}
 =\theta\Bigl(1-\frac{N_d}{n+1}\Bigr)^{d+1}-\frac{(d+1)N_d}{n+1}.
\end{equation*}
The right-hand side tends to $\theta$ as $n\to\infty$, so $c_d^{\mathrm{top}}\ge\theta$. Together with the upper bound of Bukh, Matou\v{s}ek, and Nivasch~\cite[Theorem~1.3]{BMN}, this yields $\theta\le c_d^{\mathrm{top}}\le c_d^{\mathrm{aff}}\le\theta$.
\end{proof}

\subsubsection*{Proof outline}

The proof starts with an argument due to Dotterrer, Kaufman, and Wagner~\cite{DKW} converting 
the continuous map $f$ into combinatorial data (Lemma~\ref{lem:intersection}): a triangulation $T$ of the target, compactified to a sphere, and, for every simplex $\tau$ of $T$, a mod-two cochain $a(\tau)$ on $\Delta^n$ recording which faces have images meeting~$\tau$. For a vertex $v$ of~$T$, the cochain $a(v)$ is essentially the family of $d$-faces whose images cover~$v$; it is a cocycle, and we have $|\supp a(v)|\le\ov(f)$. The cochains satisfy coboundary identities along $T$, and their sum over the top-dimensional simplices of $T$ is the constant cochain~$\one$. We show that if $\ov(f)$ is too small, then such a family of cochains cannot exist. 

The certificate of non-existence is presented in Theorem~\ref{thm:filling}. The certificate consists of a family of subsets $B_{\tau}$ of $\Delta^n$ and subcomplexes $I_{\tau}$ of $\Delta^n$ for $\tau \in T$, satisfying nesting conditions, including $|I_{\tau}| \subseteq B_{\tau}$, and such that \begin{itemize} \item the complexes $I_{\tau}$ are ``local'', in the sense that $a(\tau)$ vanishes on $I_{\tau}$; \item the sets $B_{\tau}$ are ``global'', meaning $\bigcap_{\tau \in T} B_{\tau} \neq \emptyset$; \item $B_{\tau}$ has similar topology to $|I_{\tau}|$, meaning $H_{j}(B_{\tau}, |I_{\tau}|) = 0$ for $0 \leq j \leq d$.\footnote{In fact this condition is only required for $j = d - \dim \tau$.}
\end{itemize} 

We define $I_{\tau}$ and $B_{\tau}$ by regarding $\Delta^n$ as the simplex of probability vectors on the vertex set $[n]$. For a family $G$ of $d$-faces and a probability vector $\mu$, the \emph{weighted density} $p_G(\mu)$ is the probability that $d+1$ independent samples from $\mu$ are distinct and form a member of~$G$; it is a polynomial in~$\mu$. With every simplex $\tau$ of $T$ we associate a $d$-cocycle $G_{\tau}$, which is essentially $\supp a(v)$ for a vertex $v$ of $\tau$, and (more or less) define $B_{\tau}$ as the sublevel set $\{p_{G_{\tau}} < \theta\}$ of the above polynomial and $I_{\tau}$ as the independence complex of $G_{\tau}$, so that $|I_{\tau}|= \{p_{G_{\tau}} = 0\}$. The upper bound on $\ov(f)$ needed to establish non-existence is used to show that every $B_{\tau}$ contains the barycenter of $\Delta^n$. See Section~\ref{sec:proof} for details. 

The main topological ingredient, Theorem~\ref{thm:weights}, verifies the relative homology condition above. It states that for a $d$-cocycle $G$ the relative homology of the pair $(\{p_G<\theta\},\{p_G=0\})$ vanishes in degrees at most~$d$. Its proof in Section~\ref{sec:morse} is a standard application of Morse theory, modulo a couple of technicalities. The essential input that makes the machinery applicable is a lower bound on the number of negative directions of the Hessian of $p_G$ at critical points.

Lemma~\ref{lem:index}, which provides this input, is the final novel ingredient of the proof. Its proof involves bounding eigenvalues of the Hessian after rescaling the variables and a combinatorial sublemma which shows that every positive face-critical value of $p_G$ equals $\theta$ when $G$ is a ``triangle-free'' cocycle.

The proof thus combines algebraic topology, Morse theory, linear algebra, and a touch of combinatorics. After a couple of weeks of digestion, none of the arguments any longer appear to be particularly difficult to the author, but the intuition behind some of the choices remains opaque.

\subsubsection*{Organization}

Section~\ref{sec:prelim} fixes conventions and proves the small-chain lemma. Section~\ref{sec:intersection} derives the intersection-cochain lemma from results of Dotterrer, Kaufman, and Wagner~\cite{DKW} and proves the descending filling theorem. Section~\ref{sec:weights} introduces weighted densities and deduces the overlap bound from Theorem~\ref{thm:weights}. Section~\ref{sec:morse} proves that theorem using the index bound and Morse theory. Appendix~\ref{app:example} works out a small example of the construction.

\subsubsection*{AI disclosure} This project was motivated by the high-dimensional bramble lemma~\cite{Norin}, which the author originally proved in an effort to deconstruct the application of Gromov's topological overlap~\cite{Gromov} in~\cite{EHMNSW}. A few years ago, the author observed that this lemma could be used to recover optimal overlap bounds in two dimensions~\cite{Norin}. Using a significantly more technical argument in joint work with Zden\v{e}k Dvo\v{r}\'ak, the author obtained suboptimal bounds in three dimensions, matching Gromov's.

Recently, the author tasked GPT-6 Astra with obtaining a common strengthening of the bramble and Gromov's approaches to topological overlap in three dimensions. After a few rounds of improvements, Astra was encouraged to abandon numerical optimization and prioritize elegant conceptual approaches. At this point, Astra found a proof of the optimal bound in three dimensions, which immediately extended to every dimension.

The author rewrote the proof, restructuring it and improving the presentation. Astra and Claude Opus-5.5 were enlisted to help at various points. In particular, Claude Opus-5.5 produced the figure.

\paragraph*{Lean formalization.} GPT-6.1 Sol produced a Lean formalization of Theorem~\ref{thm:main} and a strengthening of Theorem~\ref{thm:finite}, based on an earlier draft. It is available at \url{https://github.com/snorin239/SharpTopologicalOverlap}.

\section{Preliminaries}\label{sec:prelim}

In this section we introduce the standard concepts and results in algebraic topology underlying the proof and fix the notation. 

The textbook~\cite{Hatcher} by Hatcher covers all the necessary background. In fact, familiarity with a small fragment of~\cite{Hatcher} suffices. 
For standard background on homotopies, simplicial complexes, and singular and relative homology, we refer to Hatcher~\cite[Chapter~0, pp.~1--4; \S2.1, pp.~102--125 and~129]{Hatcher}. For homology with coefficients in~$\F$ and for cochains and coboundaries, see~\cite[\S2.2, p.~153; \S3.1, pp.~190--191]{Hatcher}.

\subsubsection*{Simplicial complexes}\label{sec:complexes}

Let $V$ be a finite set. A \emph{simplicial complex} on $V$ is a family $K$ of nonempty subsets of $V$, called \emph{faces}, closed under taking nonempty subsets. A \emph{$k$-face} is a face $\sigma$ of \emph{dimension} $k=|\sigma|-1$; for $k\ge1$, its \emph{facets} are the faces $\sigma\setminus\{v\}$ for $v\in\sigma$; a vertex has no facets. We identify vertices with their singleton faces, writing $v$ for $\{v\}$ when used as an index or argument. A \emph{subcomplex} is a subfamily that is a simplicial complex. The \emph{full simplex} $\Delta(V)$ consists of all nonempty subsets of~$V$.

We realize $\Delta(V)$ as the probability simplex
\begin{equation*}
 \Delta(V)=\Bigl\{\mu\in\R^V:\mu_i\ge0\ ,\ \sum_{i\in V}\mu_i=1\Bigr\},
\end{equation*}
and use the same symbol for the complex and for this set. For $i\in V$ let $e_i\in\Delta(V)$ be the corresponding unit vector. The \emph{support} of $\mu\in\R^V$ is $$\supp\mu=\{i:\mu_i\ne0\},$$ and the \emph{mass} of $A\subseteq V$ is $$\mu(A)=\sum_{i\in A}\mu_i.$$ For a nonempty $\sigma\subseteq V$, the \emph{closed face} and the \emph{open face} of $\sigma$ are
\begin{equation*}
 \Delta(\sigma)=\{\mu\in\Delta(V):\supp\mu\subseteq\sigma\},\qquad
 \Delta^\circ(\sigma)=\{\mu\in\Delta(V):\supp\mu=\sigma\};
\end{equation*}
we identify $\R^\sigma$ with the vectors in $\R^V$ that vanish outside~$\sigma$. The \emph{geometric realization} of $K$ is $|K|=\bigcup_{\sigma\in K}\Delta(\sigma)$; thus $\mu\in|K|$ if and only if $\supp\mu\in K$. The \emph{uniform vector} (or \emph{barycenter}) $b_V$ of $\Delta(V)$ has all coordinates equal to $1/|V|$.

\subsubsection*{Chains and cochains}\label{sec:chains}

All chains, cochains, and homology groups in this paper  have coefficients in~$\F$. Let $K$ be a simplicial complex. For $k\ge0$, the group of \emph{$k$-chains} $C_k(K)$ is the $\F$-vector space with basis the $k$-faces of~$K$; we identify a chain with the set of faces having coefficient one. The \emph{boundary} map $\partial:C_k(K)\to C_{k-1}(K)$ for $k \geq 1$ is
\begin{equation*}
 \partial\sigma=\sum_{v\in\sigma}\sigma\setminus\{v\},
\end{equation*}
and $\partial=0$ on $C_0(K)$. Thus our chain complexes are unaugmented, and a vertex has zero boundary. One has $\partial\partial=0$. A \emph{cycle} is a chain with zero boundary.

A \emph{$k$-cochain} is a function $c$ from the $k$-faces of $K$ to~$\F$. Such functions form the vector space $C^k(K)$. The \emph{support} $\supp c$ of $c$ is the set of $k$-faces $\sigma$ with $c(\sigma)=1$. As for chains, we often identify $c$ with $\supp c$ and call the elements of $\supp c$ the \emph{members} of~$c$; thus a $k$-cochain is a family of $k$-faces, and the sum of two cochains is their symmetric difference. A cochain is evaluated on a chain by $c(z)=\sum_{\sigma\in z}c(\sigma)$. The \emph{coboundary} $\delta:C^k(K)\to C^{k+1}(K)$ is
\begin{equation*}
 (\delta c)(\sigma)=\sum_{v\in\sigma}c(\sigma\setminus\{v\})=c(\partial\sigma);
\end{equation*}
hence $(\delta c)(z)=c(\partial z)$ for every $(k+1)$-chain $z$, and $\delta\delta=0$. A \emph{cocycle} is a cochain with $\delta c=0$. We say that a $k$-cochain $c$ \emph{vanishes on} a set $Y$ of vertices if $c(\sigma)=0$ for every $k$-face $\sigma\subseteq Y$, and that it vanishes on a simplicial subcomplex $L$ if $c(\sigma)=0$ for every $k$-face $\sigma$ of $L$. We write $\one\in C^0(K)$ for the cochain with value one on every vertex.

More generally, a \emph{chain complex} $C_*$ consists of $\F$-vector spaces $C_k$, $k\ge0$, and linear maps $\partial:C_k\to C_{k-1}$ for $k\ge1$ with $\partial\partial=0$; we put $C_k=0$ for $k<0$ and $\partial=0$ on $C_0$. An element of $C_k$ is a \emph{cycle} if its boundary is zero, and a \emph{boundary} if it lies in $\partial C_{k+1}$; every boundary is a cycle because $\partial\partial=0$. The \emph{homology group} $H_k(C_*)$ in degree~$k$ is the quotient of the subspace of cycles in $C_k$ by the subspace $\partial C_{k+1}$ of boundaries; it is an $\F$-vector space. The \emph{class} of a cycle $z\in C_k$ is its image in $H_k(C_*)$. The dual cochain complex of $C_*$ consists of the spaces $C^k=\operatorname{Hom}(C_k,\F)$ and the maps $\delta:C^k\to C^{k+1}$, $(\delta c)(z)=c(\partial z)$. A \emph{subcomplex} $D_*$ of $C_*$ consists of subspaces $D_k\subseteq C_k$ with $\partial D_k\subseteq D_{k-1}$. The \emph{quotient complex} $C_*/D_*$ consists of the quotient spaces $C_k/D_k$ and the maps $C_k/D_k\to C_{k-1}/D_{k-1}$, $z+D_k\mapsto\partial z+D_{k-1}$, which are well defined because $\partial D_k\subseteq D_{k-1}$. The \emph{relative homology} $H_k(C_*,D_*)$ is the homology in degree $k$ of the quotient complex $C_*/D_*$. An element $z\in C_k$ is a \emph{relative cycle} if $\partial z\in D_{k-1}$, that is, if $z+D_k$ is a cycle of $C_*/D_*$; the class of $z$ in $H_k(C_*,D_*)$ is the class of $z+D_k$. Thus $H_k(C_*,D_*)=0$ if and only if every $z\in C_k$ with $\partial z\in D_{k-1}$ has the form $z=\partial w+z'$ with $w\in C_{k+1}$ and $z'\in D_k$; we only use this characterization of the relative homology. 

A \emph{chain map} $\rho:C'_*\to C_*$ consists of linear maps $C'_k\to C_k$ with $\rho\partial=\partial\rho$; then $\delta(c\circ\rho)=(\delta c)\circ\rho$ for every cochain $c$ of $C_*$. A chain map $\rho$ \emph{induces} linear maps $H_k(C'_*)\to H_k(C_*)$ sending the class of $z$ to the class of $\rho(z)$. If $D'_*$ and $D_*$ are subcomplexes of $C'_*$ and $C_*$ with $\rho(D'_k)\subseteq D_k$ for every $k$, then $z+D'_k\mapsto\rho(z)+D_k$ is a chain map $C'_*/D'_*\to C_*/D_*$, and so $\rho$ induces linear maps $H_k(C'_*,D'_*)\to H_k(C_*,D_*)$. When $C'_*$ is a subcomplex of $C_*$ and $\rho$ is the inclusion, we say that these maps are \emph{induced by inclusion}. For a subcomplex $D_*$ of $C_*$, the quotient map $C_*\to C_*/D_*$ is a chain map, and we say that the maps $H_k(C_*)\to H_k(C_*,D_*)$ it induces are \emph{induced by passage to the quotient}. For a simplicial complex $K$, the chains form a chain complex $C_*(K)$, whose dual is $C^*(K)$, and $C_*(L)$ is a subcomplex of $C_*(K)$ for every simplicial subcomplex $L$ of $K$.

Let $D_*$ and $C_*$ be chain complexes with $D_i=0$ for $i>d$. Following Dotterrer, Kaufman, and Wagner~\cite[Lemma~16 and Definition~17]{DKW}, a \emph{chain--cochain map} $a:D_*\to C^{d-*}$ consists of linear maps $a:D_i\to C^{d-i}$, $0\le i\le d$, such that
\begin{equation}\label{eq:chaincochain}
 \delta a(z)=a(\partial z)\qquad \mathrm{for \; all} \quad z\in D_i,\ 0\le i\le d,
\end{equation}
since $\partial=0$ on $D_0$, this says that $\delta a(z)=0$ for $z\in D_0$. If $\rho:C'_*\to C_*$ is a chain map, then $z\mapsto a(z)\circ\rho$ is a chain--cochain map $D_*\to (C')^{d-*}$. For a simplicial complex $T$ of dimension at most $d$ and a chain--cochain map $a:C_*(T)\to C^{d-*}$, \eqref{eq:chaincochain} says that $\delta a(\tau)=\sum_{v\in\tau}a(\tau\setminus\{v\})$ for every simplex $\tau$ of $T$ with $\dim\tau\ge1$, and $\delta a(v)=0$ for every vertex $v$ of~$T$.

\subsubsection*{Singular homology}

Let $Y$ be a topological space and $k\ge0$. A \emph{singular $k$-simplex} in $Y$ is a continuous map $s:\Delta^k\to Y$. The group $S_k(Y)$ of \emph{singular $k$-chains} is the $\F$-vector space with basis the singular $k$-simplices, and $\partial s=\sum_{i=0}^ks\circ\iota_i$ for $k\ge1$, where $\iota_i:\Delta^{k-1}\to\Delta^k$ is the affine map sending $e_0,\ldots,e_{k-1}$ in order to the $e_l$ with $l\ne i$; $\partial=0$ on $S_0(Y)$. Singular chains form a chain complex~\cite[p.~108 and p.~153]{Hatcher}. 

The \emph{image} of a chain is the union of the images of the simplices with coefficient one; it is compact. For $A\subseteq Y$, $S_k(A)\subseteq S_k(Y)$ consists of the chains with image in~$A$; these subspaces form a subcomplex $S_*(A)$ of $S_*(Y)$. A continuous map $\varphi:Y\to Y'$ induces linear maps $\varphi_\#:S_k(Y)\to S_k(Y')$, $s\mapsto\varphi\circ s$, with $\varphi_\#\partial=\partial\varphi_\#$. For $y\in Y$ we write $[y]$ for the singular $0$-simplex with value~$y$.

We write $H_k(Y)=H_k(S_*(Y))$ for the singular homology of $Y$. For $A\subseteq Y$, we write $H_k(Y,A)=H_k(S_*(Y),S_*(A))$ for relative singular homology. For nonempty $Y$, the \emph{reduced homology} groups are $\widetilde H_k(Y)=H_k(Y)$ for $k\ge1$ and the quotient $\widetilde H_0(Y)$ of the space of singular $0$-chains consisting of an even number of singular $0$-simplices by $\partial S_1(Y)$; cf.~\cite[\S2.1]{Hatcher}. A subspace $A\subseteq Y$ is a \emph{strong deformation retract} of $Y$ if there is a continuous map $\Psi:Y\times[0,1]\to Y$ with $\Psi(y,0)=y$, $\Psi(Y\times\{1\})\subseteq A$, and $\Psi(a,t)=a$ for all $a\in A$ and $t\in[0,1]$.

The next lemma collects standard facts from~\cite{Hatcher} about relative singular homology.

\begin{lemma}\label{lem:prism}\leavevmode
\begin{enumerate}[label=\textup{(\alph*)}]
\item Let $Y\subseteq Y'$ be spaces and $\Psi:Y\times[0,1]\to Y'$ a continuous map with $\Psi(y,0)=y$ for $y\in Y$, and put $\Psi_1=\Psi(\cdot,1)$. There are linear maps $\Pi:S_k(Y)\to S_{k+1}(Y')$ such that, for every $c\in S_k(Y)$,
\begin{equation*}
 \partial\Pi(c)+\Pi(\partial c)=c+(\Psi_1)_\#c,
\end{equation*}
and the image of $\Pi(c)$ is contained in $\Psi(\im c\times[0,1])$.
\item If $A$ is a strong deformation retract of $Y$, then $H_k(Y,A)=0$ for all~$k$.
\item If $Y_0\subseteq Y_1\subseteq\cdots\subseteq Y_l$ and $H_k(Y_i,Y_{i-1})=0$ for $1\le i\le l$, then $H_k(Y_l,Y_0)=0$.
\end{enumerate}
\end{lemma}
\begin{proof}
Part~(a) is Hatcher's prism construction~\cite[proof of Theorem~2.10, pp.~112--113]{Hatcher}, with coefficients in~$\F$. Part~(b) follows from homotopy invariance for pairs~\cite[Proposition~2.19, p.~118]{Hatcher}, and part~(c) follows by induction from the exact sequence of a triple~\cite[p.~118]{Hatcher}.
\end{proof}

\subsubsection*{From singular to simplicial chains}\label{sec:small}

Our final preliminary result provides a tool for converting singular homology results into statements about simplicial chains.

For an open cover $\cU$ of $\Delta^n$ and $A\subseteq\Delta^n$, let $S^\cU_*(A)$ be the subcomplex of $S_*(A)$ spanned by singular simplices whose image lies in one member of $\cU$; we call these simplices \emph{small}. These chain complexes are nested under inclusions of $A$, and $S^\cU_0(A)=S_0(A)$.

We use the following result of Hatcher, which holds with coefficients in~$\F$ by the same proof as for integer coefficients~\cite[p.~153]{Hatcher}.

\begin{lemma}[{\cite[Proposition~2.21, p.~119]{Hatcher}}]\label{lem:smallhomology}
Let $Y$ be a topological space and let $\mathcal U$ be an open cover of~$Y$. Let $S_*^{\mathcal U}(Y)$ be the subcomplex of $S_*(Y)$ spanned by the singular simplices whose image lies in one member of~$\mathcal U$. Then inclusion induces isomorphisms
\begin{equation*}
 H_k(S_*^{\mathcal U}(Y))\cong H_k(Y)
\end{equation*}
for every $k \geq 0$.
\end{lemma}

\begin{lemma}\label{lem:smallchains}
Let $n\ge0$. There are an open cover $\cU$ of $\Delta^n$ and a chain map $\rho:S_*(\Delta^n)\to C_*(\Delta^n)$ such that
\begin{enumerate}[label=\textup{(\alph*)}]
\item for every $A\subseteq Y\subseteq\Delta^n$, inclusion induces isomorphisms
\begin{equation*}
 H_k(S^\cU_*(Y),S^\cU_*(A))\cong H_k(Y,A)
\end{equation*}
for each $k \geq 0$.
\item $\rho([\mu])$ is a vertex of $\Delta^n$ for every $\mu\in\Delta^n$, and $\rho(S^\cU_*(|L|))\subseteq C_*(L)$ for every subcomplex $L$ of $\Delta^n$.
\end{enumerate}
\end{lemma}

\begin{proof}
	For a nonempty $\sigma\subseteq[n]$, let $U_\sigma$ consist of all $\mu\in\Delta^n$ such that $\mu_i> 0$ for $i \in \sigma$ and $\mu_i>\mu_j$ for all $i\in\sigma,\ j\in[n]\setminus\sigma.$
	Each $U_\sigma$ is open in $\Delta^n$. If $\mu\in U_\sigma$, then $\sigma\subseteq\supp\mu$, and every coordinate at which $\mu$ attains its maximum lies in~$\sigma$. Every $\mu\in\Delta^n$ lies in $U_\sigma$ for the set $\sigma$ of coordinates at which $\mu$ attains its maximum, since this maximum is at least $1/(n+1)>0$. Thus $\cU=\{U_\sigma:\varnothing\ne\sigma\subseteq[n]\}$ is an open cover of $\Delta^n$.
	
	Fix $A\subseteq Y\subseteq\Delta^n$. Lemma~\ref{lem:smallhomology}, applied to the restricted covers of $A$ and $Y$, whose members $U_\sigma\cap A$ and $U_\sigma\cap Y$ are open and cover the respective spaces, gives homology isomorphisms for $A$ and $Y$.
	
	A short exact sequence of chain complexes induces a long exact sequence of homology groups~\cite[Theorem~2.16, p.~117]{Hatcher}; exactness means that the image of each map is the kernel of the next. For brevity, write $H_k^\cU(X)=H_k(S_*^\cU(X))$ for $X\subseteq\Delta^n$ and $H_k^\cU(Y,A)=H_k(S_*^\cU(Y),S_*^\cU(A))$. Applying this fact to the short exact sequences
	\begin{align*}
	 0&\longrightarrow S_*^\cU(A)\longrightarrow S_*^\cU(Y)\longrightarrow S_*^\cU(Y)/S_*^\cU(A)\longrightarrow0,\\
	 0&\longrightarrow S_*(A)\,\longrightarrow S_*(Y)\,\longrightarrow S_*(Y)/S_*(A)\,\longrightarrow0
	\end{align*}
	gives the two long exact sequences
	\begin{align*}
	 \cdots&\longrightarrow H_k^\cU(A)\longrightarrow H_k^\cU(Y)\longrightarrow H_k^\cU(Y,A)
	 \longrightarrow H_{k-1}^\cU(A)\longrightarrow H_{k-1}^\cU(Y)\longrightarrow\cdots,\\[3pt]
	 \cdots&\longrightarrow H_k(A)\,\longrightarrow H_k(Y)\,\longrightarrow H_k(Y,A)\,\longrightarrow H_{k-1}(A)\,\longrightarrow H_{k-1}(Y)\,\longrightarrow\cdots\,.
	\end{align*}
	The maps from absolute to absolute homology and from absolute to relative homology are induced by inclusion and passage to the quotient, respectively. The connecting map from relative homology to the homology of $A$ sends the class of a relative cycle $c$ to the class of $\partial c$. 
	
	The inclusions of small chains into singular chains give a commutative map between these short exact sequences, and hence between their long exact homology sequences by naturality~\cite[p.~127]{Hatcher}. The five lemma~\cite[p.~129]{Hatcher} states that, in a commutative diagram with two exact rows of five abelian groups, if the first, second, fourth, and fifth vertical maps are isomorphisms, then the third is an isomorphism as well. Apply it to the five consecutive terms displayed in each long exact sequence. The four maps involving the absolute homology of $A$ and $Y$ are isomorphisms by Lemma~\ref{lem:smallhomology}, so the middle map $H_k^\cU(Y,A)\to H_k(Y,A)$ is an isomorphism. (For $k=0$, the degree $-1$ groups are zero, and the same argument applies.) This proves~(a).
	
	For each $\mu\in\Delta^n$ choose a coordinate $j(\mu)$ at which $\mu$ attains its maximum. For a singular $k$-simplex $s$ with vertices $s(e_0),\ldots,s(e_k)$ define
	\begin{equation*}
		\rho(s)=\begin{cases}\{j(s(e_0)),\ldots,j(s(e_k))\}&\text{if these $k+1$ elements of $[n]$ are distinct},\\ 0&\text{otherwise},\end{cases}
	\end{equation*}
	a $k$-face of $\Delta^n$ or zero, and extend $\rho$ linearly to a map $S_*(\Delta^n)\to C_*(\Delta^n)$.
	It is routine to verify that $\rho$ is a chain map, and $\rho([\mu])=\{j(\mu)\}$ is a vertex of $\Delta^{n}$ for every $\mu \in \Delta^n$.
	
	Now let $s$ be a singular simplex with image in $|L|\cap U_\sigma$, and put $\mu=s(e_0)$. Then $\sigma\subseteq\supp\mu$, and $\supp\mu\in L$ because $\mu\in|L|$; hence $\sigma\in L$. Each $j(s(e_l))$ is a maximizing coordinate of a point of $U_\sigma$, so $j(s(e_l))\in\sigma$. Therefore $\rho(s)$ is zero or a face of $\sigma$, and it belongs to $C_*(L)$. This proves~(b).
\end{proof}

\section{Intersection cochains and descending fillings}\label{sec:intersection}

\subsection{Intersection cochains}

The following lemma  is the only input about continuous maps in the proof of Theorem~\ref{thm:finite}. It converts a continuous map $f:\Delta^n\to\R^d$ to a chain--cochain map $a:C_*(T)\to C^{d-*}(\Delta^n)$, where $T$ is a triangulation of the sphere $S^d$. The properties of this map lead to a contradiction if $\ov(f)$ is smaller than the bound prescribed by Theorem~\ref{thm:finite}. 

\begin{lemma}\label{lem:intersection}
Let $n\ge d$ and let $f:\Delta^n\to\R^d$ be continuous. There are a finite simplicial complex $T$ of dimension $d$, a $d$-cycle $\zeta\in C_d(T)$, and a chain--cochain map $a:C_*(T)\to C^{d-*}(\Delta^n)$ with the following properties.
\begin{enumerate}[label=\textup{(\alph*)}]
\item\label{I:top} $a(\zeta)=\one$ in $C^0(\Delta^n)$.
\item\label{I:vertex} For every vertex $v$ of $T$, $|\supp a(v)|\le\ov(f)$.
\item\label{I:gp} For every $\tau\in T$ of dimension $i\ge1$, $\supp a(\tau)$ does not contain $\lfloor d/i\rfloor+1$ pairwise vertex-disjoint faces.
\end{enumerate}
\end{lemma}

The complex $T$, the cycle $\zeta$, and the map $a$ are constructed by Dotterrer, Kaufman, and Wagner~\cite{DKW}; the properties we need follow from facts proved there, which are spread over several lemmas, and we give pointers to the relevant passages below. In their construction, $f$ is approximated by a piecewise linear map in general position, $T$ is a fine triangulation of the sphere $S^d=\R^d\cup\{\infty\}$, and for a simplex $\tau$ of $T$ and a face $\sigma$ of $\Delta^n$ of complementary dimension, $a(\tau)(\sigma)$ is the parity of the number of points of $\Delta(\sigma)$ that the approximation maps into~$\tau$. Cochains of this kind were also used by Gromov~\cite[\S2.2]{Gromov}.

\begin{proof}[Proof of Lemma~\ref{lem:intersection}]
We apply the construction at the beginning of the proof of~\cite[Theorem~8, p.~314]{DKW} to the restriction of $f$ to the $d$-skeleton of $\Delta^n$, regarded as a map into $S^d=\R^d\cup\{\infty\}$ (cf.~\cite[Remark~9]{DKW}); $a$ is the intersection number homomorphism~\cite[(8), p.~313]{DKW}. We choose the approximation so close to $f$ that, under it, no point lies in the images of more than $\ov(f)$ $d$-faces~\cite[Lemma~10]{DKW}, and so that it misses a vertex $\infty$ of $T$; then $a(\infty)=0$ (cf.~\cite[p.~314]{DKW}). The sum $\zeta$ of the $d$-simplices of $T$ is a cycle with $a(\zeta)=\one$~\cite[p.~315]{DKW}. By~\cite[Lemma~16]{DKW}, \eqref{eq:chaincochain} holds for simplices of positive dimension. For a vertex $v$ of $T$, a path $c$ in $T$ from $\infty$ to $v$ gives $a(v)=\delta a(c)$, as in~\cite[p.~315]{DKW}, so $\delta a(v)=0$ in $C^{d+1}(\Delta^n)$. Every member of $\supp a(v)$ is a $d$-face whose image under the approximation contains $v$, which gives~\ref{I:vertex}. Property~\ref{I:gp} is shown in the proof of~\cite[Lemma~14, p.~313]{DKW}: the images of the members of $\supp a(\tau)$ have a common point, so by strong general position at most $d/i$ of them are pairwise vertex-disjoint.
\end{proof}

The reader may find it helpful to read Appendix~\ref{app:example-cochains} at this point; it describes the complex $T$ and the cochains $a(\tau)$ for an affine map $f:\Delta^4\to\R^2$, which serves as a running example.

\subsection{Descending fillings}\label{sec:filling}

Our next theorem describes a certificate of non-existence of a map $a:C_*(T)\to C^{d-*}(\Delta^n)$ satisfying the conditions of Lemma~\ref{lem:intersection}. It substantially strengthens a high-dimensional bramble lemma of the author~\cite[Theorem 1]{Norin}, but retains essentially the same proof. 

In the proof below we fill chains in the source simplex from the $d$-dimensional target simplices downwards, using relative-homology vanishing. For comparison, Gromov's bound is obtained via a cofilling argument~\cite{Gromov} which proceeds in the opposite direction, constructing cochains on the source simplex from the target vertices upwards and controlling their supports by cofilling inequalities; see also Matou\v{s}ek and Wagner~\cite[\S3.4]{MW} and Dotterrer, Kaufman, and Wagner~\cite[proof of Theorem~8]{DKW}.

\begin{theorem}\label{thm:filling}
Let $n\ge0$, let $T$ be a finite simplicial complex of dimension at most $d$, and let $a:C_*(T)\to C^{d-*}(\Delta^n)$ be a chain--cochain map. For every $\tau\in T$ let $B_\tau\subseteq\Delta^n$ be a subset and let $I_\tau$ be a simplicial subcomplex of $\Delta^n$ with $|I_\tau|\subseteq B_\tau$. Suppose that
\begin{enumerate}[label=\textup{(\roman*)}]
\item $B_\tau\subseteq B_\eta$ and $I_\tau\subseteq I_\eta$ for every face $\eta$ of $\tau$,
\item $a(\tau)$ vanishes on $I_\tau$, and
\item $H_{d-\dim\tau}(B_\tau,|I_\tau|)=0$.
\end{enumerate}
If $\bigcap_{\tau\in T}B_\tau\ne\varnothing$, then $a(\zeta)=0$ for every $d$-cycle $\zeta\in C_d(T)$.
\end{theorem}

\begin{proof}
Let $\zeta\in C_d(T)$ be a cycle, choose $b\in\bigcap_{\tau\in T}B_\tau$, and let an open cover $\cU$ and $\rho:S_*(\Delta^n)\to C_*(\Delta^n)$  be given by Lemma~\ref{lem:smallchains}. By~(iii) and Lemma~\ref{lem:smallchains}(a), for $\tau \in T$ we have
\begin{equation*}
 H_{d-\dim\tau}(S^\cU_*(B_\tau),S^\cU_*(|I_\tau|))=0.
\end{equation*}
For $i=d,d-1,\ldots,0$ in turn and every $i$-simplex $\tau$ of $T$ we choose singular chains
\begin{equation*}
 z_\tau\in S^\cU_{d-i}(B_\tau),\qquad
 c_\tau\in S^\cU_{d-i+1}(B_\tau),\qquad
 z_\tau-\partial c_\tau\in S^\cU_{d-i}(|I_\tau|),
\end{equation*}
as follows. For a $d$-simplex $\tau$, let $z_\tau=[b]$ if $\tau\in\zeta$ and $z_\tau=0$ otherwise. For an $i$-simplex $\eta$ with $i<d$, put
\begin{equation*}
 z_\eta=\sum_{\tau\supset\eta,\ \dim\tau=i+1}c_\tau;
\end{equation*}
it lies in $S^\cU_{d-i}(B_\eta)$ by~(i). In each case we show $\partial z_\eta\in S^\cU_{d-i-1}(|I_\eta|)$; the relative-homology vanishing above then provides $c_\eta$. For $i=d$, $z_\eta$ is a $0$-chain and $\partial z_\eta=0$. For $i<d$,
\begin{equation*}
 \partial z_\eta=\sum_{\tau\supset\eta,\ \dim\tau=i+1}(\partial c_\tau-z_\tau)+\sum_{\tau\supset\eta,\ \dim\tau=i+1}z_\tau.
\end{equation*}
The first sum lies in $S^\cU_{d-i-1}(|I_\eta|)$ by~(i). If $i=d-1$, the second sum is $[b]$ times the number of $d$-simplices of $\zeta$ containing $\eta$, which is even because $\zeta$ is a cycle. If $i\le d-2$, the second sum is the sum of $c_\xi$ over the pairs $(\tau,\xi)$ of simplices of $T$ with $\eta\subset\tau\subset\xi$, $\dim\tau=i+1$, and $\dim\xi=i+2$; every $\xi$ occurs in exactly two pairs, since $|\xi\setminus\eta|=2$. In both cases the second sum vanishes.

For $0\le i\le d$ let
\begin{equation*}
 E_i=\sum_{\dim\tau=i}a(\tau)(\rho(z_\tau))\in\F.
\end{equation*}
Let $\tau$ be an $i$-simplex of $T$. Lemma~\ref{lem:smallchains}(b) gives $\rho(z_\tau-\partial c_\tau)\in C_{d-i}(I_\tau)$, so~(ii) gives $a(\tau)(\rho(z_\tau-\partial c_\tau))=0$. Since $\rho$ is a chain map, \eqref{eq:chaincochain} yields
\begin{align*}
 a(\tau)(\rho(z_\tau))
 &=a(\tau)(\rho(\partial c_\tau))
 =a(\tau)(\partial\rho(c_\tau))\\
 &=(\delta a(\tau))(\rho(c_\tau))
 =a(\partial\tau)(\rho(c_\tau)).
\end{align*}
For $i\ge1$ this gives
\begin{align*}
 E_i
 &=\sum_{\dim\tau=i}\ \sum_{v\in\tau}a(\tau\setminus\{v\})(\rho(c_\tau))\\
 &=\sum_{\dim\eta=i-1}a(\eta)\Bigl(\rho\Bigl(\sum_{\tau\supset\eta,\ \dim\tau=i}c_\tau\Bigr)\Bigr)
 =E_{i-1}.
\end{align*}
For $i=0$, $\delta a(v)=0$ for every vertex $v$, so $E_0=0$. Thus
\begin{equation*}
 a(\zeta)(\rho([b]))=\sum_{\tau\in\zeta}a(\tau)(\rho([b]))=E_d=E_{d-1}=\cdots=E_0=0.
\end{equation*}
By Lemma~\ref{lem:smallchains}(b), $\rho([b])=\{x\}$ for some $x\in[n]$. Also $\delta a(\zeta)=a(\partial\zeta)=0$, so for distinct $u,v\in[n]$,
\begin{equation*}
 a(\zeta)(u)+a(\zeta)(v)=(\delta a(\zeta))(\{u,v\})=0.
\end{equation*}
Hence $a(\zeta)$ is constant on $[n]$, and $a(\zeta)(x)=0$ implies $a(\zeta)=0$.
\end{proof}

\section{Weighted densities and the overlap bound}\label{sec:weights}\label{sec:proof}

In this section, given a complex $T$ and a chain--cochain map $a$ as in Lemma~\ref{lem:intersection}, coming from a continuous map $f: \Delta^n \to \R^d$ with $\ov(f)$ smaller than the bound in Theorem~\ref{thm:finite}, we construct sets $I_{\tau}, B_{\tau}$ as in Theorem~\ref{thm:filling}, modulo Theorem~\ref{thm:weights}, which is proved in Section~\ref{sec:morse}, and provides the relative homology condition in Theorem~\ref{thm:filling}(iii). The resulting contradiction implies Theorem~\ref{thm:finite}.

\subsection{The main topological theorem}\label{sec:polynomials}

The sets $B_\tau$ are defined using weighted density functions associated with families of $d$-faces of $\Delta^n$. We now describe these functions.

Let $V$ be a nonempty finite set. For a family $G$ of $(d+1)$-element subsets of $V$ and $\mu=(\mu_i)_{i\in V}\in\R^V$, define
\begin{equation}\label{eq:polynomial}
 p_G(\mu)=(d+1)!\sum_{e\in G}\prod_{i\in e}\mu_i.
\end{equation}
The polynomial $p_G$ is multilinear and homogeneous of degree $d+1$. It has nonnegative coefficients, so it is nondecreasing in each coordinate on $[0,\infty)^V$.  For $\mu\in\Delta(V)$, $p_G(\mu)$ is the probability that $d+1$ independent samples from $\mu$ are pairwise distinct and their set belongs to $G$. We call $p_G(\mu)$ the \emph{weighted density} of $G$ at~$\mu$. In particular, $0\le p_G\le1$ on $\Delta(V)$, and at the uniform vector $b_V$,
\begin{equation}\label{eq:uniform}
 p_G(b_V)=\frac{(d+1)!\,|G|}{|V|^{d+1}}.
\end{equation}
%If $|V|\ge d+1$, then $p_G(b_V)\le |G|/\binom{|V|}{d+1}$.

We have $p_G(\mu)=0$ if and only if $\supp\mu$ contains no member of $G$, that is,
\begin{equation}\label{eq:independence}
 \{\mu\in\Delta(V):p_G(\mu)=0\}=|\Ind(G)|,
\end{equation}
where the \emph{independence complex} $\Ind(G)$ of $G$ consists of the nonempty sets $\sigma\subseteq V$ that contain no member of~$G$. 

For a function $\varphi$ on $\Delta(V)$ we write $\{\varphi<t\}=\{\mu\in\Delta(V):\varphi(\mu)<t\}$, and similarly $\{\varphi\le t\}$ and $\{\varphi=0\}$. The following is the main novel topological ingredient of our proof.

\begin{theorem}\label{thm:weights}
Let $V$ be a nonempty finite set and let $G$ be a $d$-cocycle of $\Delta(V)$. Then
\begin{equation*}
 H_j\bigl(\{p_G<\theta\},\{p_G=0\}\bigr)=0
\end{equation*}
for $0 \leq j \leq d$.
\end{theorem}

%By the characterization of relative homology in ``\nameref{sec:chains}'', Theorem~\ref{thm:weights}  implies that every singular chain of dimension $j\le d$ in $\{p_G<\theta\}$ whose boundary lies in $|\Ind(G)|=\{p_G=0\}$ is homologous, inside $\{p_G<\theta\}$, to a chain in $|\Ind(G)|$.

Note that the threshold $\theta$ in Theorem~\ref{thm:weights} cannot be replaced by a larger number and the inequality cannot be made non-strict. Consider for example $V=\{1,\ldots,d+1\}$ and $G=\{V\}$. Then $\{p_G=0\}$ is the boundary of the $d$-simplex $\Delta(V)$, and $\{p_G<\theta\}$ is $\Delta(V)$ with its uniform vector removed, so Theorem~\ref{thm:weights} holds in this case  by Lemma~\ref{lem:prism}(b). On the other hand, $\{p_G\le\theta\}=\Delta(V)$, and $H_d(\Delta(V),\partial\Delta(V))\cong \widetilde H_{d-1}(\partial\Delta(V))\cong\F$.

The proof of  Theorem~\ref{thm:weights} is postponed until Section~\ref{sec:morse}.

\subsection{Proof of Theorem~\ref{thm:finite}.}

We continue working towards the proof of Theorem~\ref{thm:finite}.  
Throughout this section let $n+1>N_d$, let $f:\Delta^n\to\R^d$ be a continuous map, put $m=\ov(f)$, and let $T$, $\zeta$, and $a$ be given by Lemma~\ref{lem:intersection} applied to $f$. (Since $N_d\ge d$ and $n+1>N_d$, we have $n\ge d$, so Lemma~\ref{lem:intersection} applies.)

Our immediate goal is defining regions $B_\tau$ and subcomplexes $I_\tau$ of $\Delta^n$ for each $\tau \in T$ satisfying the conditions of Theorem~\ref{thm:filling}. The pairs $(B_\tau,|I_\tau|)$ will essentially (but not exactly) correspond to the pairs considered in Theorem~\ref{thm:weights}.

\subsubsection*{Outline of the construction.}
We outline a naive construction first, and then explain how our final choice repairs the issues that arise.  

The most straightforward way of defining $B_{\tau}$ and $I_{\tau}$ based on Theorem~\ref{thm:weights} would be to choose for every $\tau \in T$ a family $G_{\tau}$ of $(d+1)$-element subsets of $[n]$ and set $$B_{\tau} = \{\mu \in \Delta^n \,|\, p_{G_{\tau}}(\mu) < \theta\} \quad \mathrm{and} \quad I_{\tau} = \Ind(G_{\tau}).$$
If $G_{\tau}$ is a cocycle then Theorem~\ref{thm:weights} guarantees condition (iii) of Theorem~\ref{thm:filling}.

The choice of $G_{\tau}$ is motivated by the remaining conditions of Theorem~\ref{thm:filling}.
The nesting requirement $B_{\tau} \subseteq B_{\eta}$ for $\eta \subseteq \tau$ in Theorem~\ref{thm:filling}(i) suggests choosing $G_{\eta} \subseteq G_{\tau}$ for $\eta \subseteq \tau$ to ensure $p_{G_{\eta}}(\mu) \leq p_{G_{\tau}}(\mu)$ for all $\mu \in \Delta^n$. Thus, in particular, we would like $G_{\tau} \supseteq \cup_{v \in \tau} G_{v}$ and choosing 
$$G_{\tau} = \cup_{v \in \tau} G_{v}$$ guarantees the nesting and reduces the choice of $G_{\tau}$ to the choice at vertices.

Theorem~\ref{thm:filling}(ii) requires $a(v)$ to vanish on $I_{v} =  \Ind(G_{v})$, i.e. $G_{v}$ needs to contain $\supp(a(v))$
and defining $$G_v = \supp(a(v))$$ would complete our initial naive choice. Note that $a(v)$ is conveniently a cocycle as $a$ is a chain--cochain map. 

The immediate issues with this choice are the following: \begin{itemize}
	\item $\Ind(G_\tau)$ contains the complete $(d-1)$-skeleton of $\Delta^n$, so when $\dim\tau\ge1$, a nonzero cochain $a(\tau)$ cannot vanish on it;
	\item the union of cocycles $G_{\tau}$ need not be a cocycle.
\end{itemize}

We can repair both of these issues by restricting every $G_{\tau}$ (and thus also $B_{\tau}$ and $I_{\tau}$)  to a vertex set $V_{\tau}$ on which all cochains $a(\eta)$ vanish for $\eta \subseteq \tau$ with the dimension of $\eta$ at least one. Indeed, the chain--cochain conditions ensure that the cochains $a(v)$ for $v \in \tau$  then all agree on $\Delta(V_{\tau})$. Moreover, Lemma~\ref{lem:intersection}\ref{I:gp} ensures that $V_{\tau}$ can be obtained from $[n]$ by deleting a set of vertices of size at most $N_d$ which will be essential in our argument. 

The repair causes a new problem, not present before: after restriction the condition $\cap_{\tau \in T}B_{\tau} \neq \emptyset$ need not be satisfied.
We instead define $B_\tau$ in all of $\Delta^n$ by a homogeneous weighted density inequality that preserves the required nesting. Normalizing the restriction to $V_\tau$ gives a strong deformation retraction of $B_\tau$ onto $\{p_{G_\tau}<\theta\}\subseteq\Delta(V_\tau)$, where we apply Theorem~\ref{thm:weights}. This finishes the outline and we proceed to the precise description of the construction. 
 
\vskip 10pt 
 
Let $F_\tau$ be the set of $d$-faces of $\Delta^n$ that contain an element of $\supp a(\eta)$ for some face $\eta$ of $\tau$.
If $\eta \subseteq \tau$, then $F_\eta\subseteq F_\tau$.

For every simplex $\eta$ of $T$ of positive dimension, fix for the rest of the proof a maximal family of pairwise vertex-disjoint elements of $\supp a(\eta)$. Lemma~\ref{lem:intersection}\ref{I:gp} ensures that this family has size at most $\lfloor d/i \rfloor$, where $i$ is the dimension of $
\eta$. For $\tau\in T$, let $S_\tau\subseteq[n]$ be the set of vertices of the faces in the families chosen for the faces of $\tau$. 

Define $$V_\tau=[n]\setminus S_\tau \qquad \mathrm{and} \qquad  G_\tau=\{\sigma\in\supp a(v):\sigma\subseteq V_\tau\},$$ where $v$ is any vertex of $\tau$. The next lemma summarizes the consequences of our choices, and in particular shows that $G_{\tau}$ does not depend on the choice of $v$, and bounds $|F_\tau|$.

\begin{lemma}\label{lem:covers}
Let $\tau\in T$.
\begin{enumerate}[label=\textup{(\roman*)}]
\item $|S_\tau|\le N_d$ and $V_\tau\ne\varnothing$.
\item $S_\eta \subseteq S_\tau$ for every face $\eta$ of~$\tau$.
\item For every face $\eta$ of $\tau$ of positive dimension, $a(\eta)$ vanishes on $V_\tau$.
\item The cochains $a(v)$, $v\in\tau$, agree on the $d$-faces contained in $V_\tau$, $G_\tau$ does not depend on the vertex $v$ of $\tau$ and $G_{\tau}$ is a $d$-cocycle of~$\Delta(V_\tau)$.
\item The members of $F_\tau$ contained in $V_\tau$ are exactly the members of~$G_\tau$.
\item $|F_\tau|\le m+N_d\binom{n}{d}$.
\end{enumerate}
\end{lemma}
\begin{proof}
Let $\eta$ be a face of $\tau$ of dimension $i\ge1$. The elements of $\supp a(\eta)$ have dimension $d-i$. By Lemma~\ref{lem:intersection}\ref{I:gp}, the family chosen for $\eta$ has at most $\lfloor d/i\rfloor$ members, which have at most $(d-i+1)\lfloor d/i\rfloor$ vertices in total. By maximality, every element of $\supp a(\eta)$ shares a vertex with a member of the family, hence meets $S_\tau$; this proves~(iii). Since $\dim\tau\le d$, the simplex $\tau$ has at most $\binom{d+1}{i+1}$ faces of dimension $i$, so $$|S_\tau|\le  \sum_{i=1}^{d}\binom{d+1}{i+1}(d+1-i)\left\lfloor\frac di\right\rfloor = N_d < n+1 $$ by~\eqref{eq:C}, which proves~(i). Property~(ii) holds because every face of a face of $\tau$ is a face of~$\tau$.

 Let $\{v,w\}\subseteq\tau$ be an edge. By~\eqref{eq:chaincochain}, $a(v)+a(w)=a(\partial\{v,w\})=\delta a(\{v,w\})$, and for a $d$-face $\sigma\subseteq V_\tau$ we have $(\delta a(\{v,w\}))(\sigma)=\sum_{i\in\sigma}a(\{v,w\})(\sigma\setminus\{i\})=0$ by~(iii). This proves the first two assertions. If $W\subseteq V_\tau$ has $d+2$ elements, \eqref{eq:chaincochain} gives $\sum_{i\in W}a(v)(W\setminus\{i\})=(\delta a(v))(W)=0$, so $W$ contains an even number of members of~$G_\tau$. This proves (iv).

Finally, let $\sigma \subseteq V_\tau$ be a $d$-face. By~(iii), $\sigma$ contains no element of $\supp a(\eta)$ for a face $\eta$ of $\tau$ of positive dimension. Hence $\sigma\in F_\tau$ if and only if $a(v)(\sigma)=1$ for some $v\in\tau$. By~(iv), this means $\sigma \in G_\tau$. This proves (v).

To prove~(vi), let $v$ be a vertex of $\tau$. By~(v), the members of $F_\tau$ contained in $V_\tau$ are the members of $G_\tau\subseteq\supp a(v)$, so by Lemma~\ref{lem:intersection}\ref{I:vertex} there are at most $m$ of them. Every other member of $F_\tau$ contains a vertex of $S_\tau$. Each vertex of $\Delta^n$ lies in $\binom{n}{d}$ $d$-faces, and $|S_\tau|\le N_d$ by~(i), so at most $N_d\binom nd$ members of $F_\tau$ meet~$S_\tau$. This proves (vi).
\end{proof}

For $\tau\in T$ define
\begin{equation}\label{eq:weightspace}
 B_\tau=\bigl\{\mu\in\Delta^n:p_{F_\tau}(\mu)<\theta\,\mu(V_\tau)^{d+1}\bigr\},\qquad
 I_\tau=\Ind(G_\tau).
\end{equation}
Here by convention $p_{F_\tau}$ is the weighted density~\eqref{eq:polynomial} of $F_\tau$ on $V=[n]$, and $\Ind(G_\tau)$ is the independence complex~\eqref{eq:independence} of $G_\tau$ on $V_\tau$, a subcomplex of $\Delta^n$ whose faces are contained in~$V_\tau$. For $\mu\in\Delta(V_\tau)$, only the members of $F_\tau$ contained in $V_\tau$ contribute to $p_{F_\tau}(\mu)$, so Lemma~\ref{lem:covers}(v) gives
\begin{equation}\label{eq:restrict}
 p_{F_\tau}(\mu)=p_{G_\tau}(\mu) \quad \mathrm{for} \quad \mu \in \Delta(V_\tau).
\end{equation}

The reader may find it helpful to read Appendix~\ref{app:example-regions} at this point; it describes the sets $S_\tau$, $G_\tau$, $I_\tau$, and $B_\tau$ for the running example.

The next lemma establishes properties (i) and (iii) of Theorem~\ref{thm:filling} and the required common-point condition.

\begin{lemma}\label{lem:carriers}
Let $\tau\in T$.
\begin{enumerate}[label=\textup{(\roman*)}]
\item $|I_\tau|\subseteq B_\tau$; moreover $B_\tau\subseteq B_\eta$ and $I_\tau\subseteq I_\eta$ for every face $\eta$ of~$\tau$.
\item If $m<((n+1-N_d)/(d+1))^{d+1}-N_d\binom{n}{d}$, then the uniform vector $b_{[n]}$ of $\Delta^n$ belongs to $B_\tau$.
\item $H_j(B_\tau,|I_\tau|)=0$ for $0\le j\le d$.
\end{enumerate}
\end{lemma}
\begin{proof}
(i) A point $\mu\in |I_\tau|$ is supported in $V_\tau$, so $\mu(V_\tau)=1$ and $p_{F_\tau}(\mu)=p_{G_\tau}(\mu)=0<\theta$ by~\eqref{eq:restrict}. Let $\eta$ be a face of $\tau$. By Lemma~\ref{lem:covers}(ii), $V_\tau\subseteq V_\eta$, and since $F_\eta\subseteq F_\tau$, every $\mu\in B_\tau$ satisfies
\begin{equation*}
 p_{F_\eta}(\mu)\le p_{F_\tau}(\mu)<\theta\,\mu(V_\tau)^{d+1}\le\theta\,\mu(V_\eta)^{d+1}.
\end{equation*}
Let $\sigma\in I_\tau$. Then $\sigma\subseteq V_\tau\subseteq V_\eta$. If some $e\in G_\eta$ were contained in $\sigma$, then $e\in F_\eta\subseteq F_\tau$ and $e\subseteq V_\tau$, so Lemma~\ref{lem:covers}(v) would give $e\in G_\tau$, contrary to $\sigma\in I_\tau$. Thus $\sigma\in I_\eta$.

(ii) By~\eqref{eq:uniform} and Lemma~\ref{lem:covers}(vi) and~(i),
\begin{align*}
 p_{F_\tau}(b_{[n]})&=\frac{(d+1)!\,|F_\tau|}{(n+1)^{d+1}}\le\frac{(d+1)!\,\bigl(m+N_d\binom nd\bigr)}{(n+1)^{d+1}}\\
 &<\frac{(d+1)!}{(d+1)^{d+1}}\left(\frac{n+1-N_d}{n+1}\right)^{d+1}\le\theta\left(\frac{|V_\tau|}{n+1}\right)^{d+1}=\theta\,b_{[n]}(V_\tau)^{d+1}.
\end{align*}

(iii) We first deform $B_\tau$ onto $B_\tau\cap\Delta(V_\tau)$, keeping this intersection fixed. For $\mu\in B_\tau$ put $t=\mu(S_\tau)=1-\mu(V_\tau)$; then $t<1$, since for $t=1$ the right-hand side in~\eqref{eq:weightspace} vanishes while $p_{F_\tau}\ge0$.
For $Y\subseteq[n]$, let $\mu|_Y$ agree with $\mu$ on $Y$ and vanish elsewhere.
For $\mu\in B_\tau,\ r\in[0,1]$ define 
\begin{equation*}
 \Psi(\mu,r)=\frac{1-(1-r)t}{1-t}\,\mu|_{V_\tau}+(1-r)\,\mu|_{S_\tau}\qquad
\end{equation*}
This map is continuous, $\Psi(\mu,0)=\mu$, $\Psi(\mu,1)\in\Delta(V_\tau)$, and $\Psi(\mu,r)=\mu$ for $r \in [0,1]$ if $\mu\in B_\tau\cap\Delta(V_\tau)$.

We show $\Psi(\mu,r)\in B_\tau$. Put $s=(1-r)t\le t$; then $\Psi(\mu,r)(S_\tau)=s$. If $t=0$, then $\Psi(\mu,r)=\mu$. If $t>0$, let $\nu=\mu|_{V_\tau}/(1-t)$ and $\omega=\mu|_{S_\tau}/t$, so that $\mu=(1-t)\nu+t\omega$ and $\Psi(\mu,r)=(1-s)\nu+s\omega$. Since $p_{F_\tau}$ is homogeneous of degree $d+1$,
\begin{equation*}
 \frac{p_{F_\tau}\bigl((1-s)\nu+s\omega\bigr)}{(1-s)^{d+1}}=p_{F_\tau}\left(\nu+\frac s{1-s}\,\omega\right)\qquad(0\le s<1),
\end{equation*}
which is nondecreasing in $s$, because $s/(1-s)$ is increasing and $p_{F_\tau}$ is nondecreasing in each coordinate on $[0,\infty)^{n+1}$. Hence $p_{F_\tau}(\Psi(\mu,r))/(1-s)^{d+1}\le p_{F_\tau}(\mu)/(1-t)^{d+1}<\theta$, that is, $\Psi(\mu,r)\in B_\tau$. Thus $B_\tau\cap\Delta(V_\tau)$ is a strong deformation retract of $B_\tau$, and Lemma~\ref{lem:prism}(b) gives $$H_j(B_\tau,B_\tau\cap\Delta(V_\tau))=0$$ for all~$j$.

By~\eqref{eq:restrict}, with $\Delta(V_\tau)$ regarded as the probability simplex on $V_\tau$,
\begin{equation*}
 B_\tau\cap\Delta(V_\tau)=\{p_{G_\tau}<\theta\},\qquad |I_\tau|=\{p_{G_\tau}=0\}.
\end{equation*}
By Lemma~\ref{lem:covers}(i) and~(iv), Theorem~\ref{thm:weights} applies to $V=V_\tau$ and $G=G_\tau$, and gives $H_j(B_\tau\cap\Delta(V_\tau),|I_\tau|)=0$ for $0\le j\le d$. Lemma~\ref{lem:prism}(c), applied to $|I_\tau|\subseteq B_\tau\cap\Delta(V_\tau)\subseteq B_\tau$, proves~(iii).
\end{proof}

We are now ready to finish the proof of Theorem~\ref{thm:finite}.

\begin{proof}[Proof of Theorem~\ref{thm:finite}]
Suppose that $m<((n+1-N_d)/(d+1))^{d+1}-N_d\binom{n}{d}$. We apply Theorem~\ref{thm:filling} to $a$ and the sets $B_\tau$ and subcomplexes $I_\tau$ defined in~\eqref{eq:weightspace}.

Theorem~\ref{thm:filling}(i) holds by Lemma~\ref{lem:carriers}(i),
Theorem~\ref{thm:filling}(ii) holds by Lemma~\ref{lem:covers}(iii) and (iv), and Theorem~\ref{thm:filling}(iii) holds by Lemma~\ref{lem:carriers}(iii).

By Lemma~\ref{lem:carriers}(ii), $\cap_{\tau \in T} B_\tau \neq \emptyset$. Thus all the conditions of Theorem~\ref{thm:filling} are satisfied implying $a(\zeta)=0$, contradicting Lemma~\ref{lem:intersection}\ref{I:top}. 

Hence $\ov(f)=m\ge((n+1-N_d)/(d+1))^{d+1}-N_d\binom{n}{d}$.
\end{proof}

\section{Topology of the weighted sublevel sets}\label{sec:morse}

In this section we prove Theorem~\ref{thm:weights}.

\subsection{Critical points and the index bound}\label{sec:calculus}\label{sec:spectral}

As we describe in the next  subsection, application of Morse theory in Theorem~\ref{thm:weights} requires us to bound the negative index of $p_G$ from below at its face-critical points. This is accomplished in this subsection.

We start by defining face-critical points and the negative index.

Let $V$ be a nonempty finite set. For a nonempty $\sigma\subseteq V$, put
\begin{equation*}
 \R^\sigma_0=\Bigl\{x\in\R^V:x_i=0\ (i\notin\sigma),\ \sum_ix_i=0\Bigr\}.
\end{equation*}
The open face $\Delta^\circ(\sigma)= \{\mu\in\Delta(V):\supp\mu=\sigma\}$ is an open subset of the affine subspace $\mu+\R^\sigma_0$ for any $\mu\in\Delta^\circ(\sigma)$ with tangent space $\R^\sigma_0$ at every point. Let $\varphi$ be a $C^2$ function on an open subset of $\R^V$, with differential $d\varphi(\mu)$ and Hessian form $\Hess\varphi(\mu) = (\partial_i \partial_j \varphi (\mu))_{i,j \in V}$. A point $\mu\in\Delta(V)$ in the domain of $\varphi$, with support $\sigma$, is \emph{face-critical} for $\varphi$ if $d\varphi(\mu)$ vanishes on $\R^\sigma_0$, that is, if $\mu$ is a critical point of the restriction of $\varphi$ to $\Delta^\circ(\sigma)$. The \emph{negative index} of $\varphi$ at a face-critical point $\mu$ with support $\sigma$ is the largest dimension of a subspace of $\R^\sigma_0$ on which $\Hess\varphi(\mu)$ is negative definite; equivalently, the number of negative eigenvalues, counted with multiplicity, of $\Hess\varphi(\mu)$ restricted to $\R^\sigma_0$.  %A face-critical point is \emph{nondegenerate} if the restricted Hessian is nondegenerate.

In this section we investigate the face-critical points of the polynomial $p_G$ for $G\subseteq\binom{V}{d+1}$ defined  in Section~\ref{sec:polynomials}.
For a family $G\subseteq\binom{V}{d+1}$, $i\in V$, and $\mu\in\Delta(V)$, we have
\begin{equation}\label{eq:derivbound}
 \partial_ip_{G}(\mu)=(d+1)\cdot d!\sum_{e\in G,\ e\ni i}\ \prod_{j\in e\setminus\{i\}}\mu_j\;\in\;[0,d+1],
\end{equation}
where the upper bound holds since $d!$ times the sum is the probability that $d$ independent samples from $\mu$ are pairwise distinct and their set has the form $e\setminus\{i\}$ for some $e\in G$ containing $i$.

As $p_G$ is multilinear, the diagonal entries $\partial_i^2p_G$ of the Hessian vanish, and for distinct $i,j\in V$ and $\mu\in\Delta(V)$
\begin{equation}\label{eq:secondderiv}
 \partial_i\partial_jp_G(\mu)=(d+1)!\sum_{e\in G,\ e\supseteq\{i,j\}}\ \prod_{l\in e\setminus\{i,j\}}\mu_l.
\end{equation}

Before proceeding to the main result of this subsection, we need a couple of lemmas about face-critical points of the polynomial $p_G$.

The first one is a standard consequence of face-criticality for multilinear polynomials.

\begin{lemma}\label{lem:first}
	Let  $G\subseteq\binom{V}{d+1}$ and let $\mu\in\Delta(V)$ with $p_G(\mu) > 0$ such that $\supp \mu = V$ be face-critical for $p_G$. Let $A=(\partial_i\partial_jp_G(\mu))_{i,j\in V}$. Then, for every $i\in V$,
	\begin{enumerate}[label=\textup{(\alph*)}]
		\item $\partial_ip_G(\mu)=(d+1)p_G(\mu)$, and
		\item $\sum_{j\in V}A_{ij}\mu_j=d(d+1)p_G(\mu)$.
	\end{enumerate}
\end{lemma}

\begin{proof}
	(a) For distinct $i,j\in V$, the vector $e_i-e_j$ lies in $\R^V_0$, so face-criticality gives $\partial_ip_G(\mu)=\partial_jp_G(\mu)$. Thus all partial derivatives of $p_G$ at $\mu$ have a common value~$c$ and
	\begin{equation*}
		c=\sum_{i\in V}\mu_i\,\partial_ip_G(\mu)=(d+1)!\sum_{e\in G}(d+1)\prod_{l\in e}\mu_l=(d+1)p_G(\mu),
	\end{equation*}
	where the first equality uses $\sum_i\mu_i=1$.
	
	(b) Fix $i\in V$. Since $A_{ii}=0$, \eqref{eq:secondderiv} gives
	\begin{align*}
		\sum_{j\in V}A_{ij}\mu_j&=(d+1)!\sum_{j\ne i}\ \sum_{e\in G,\ e\supseteq\{i,j\}}\mu_j\prod_{l\in e\setminus\{i,j\}}\mu_l\\ &=(d+1)!\sum_{e\in G,\ e\ni i}\ \sum_{j\in e\setminus\{i\}}\ \prod_{l\in e\setminus\{i\}}\mu_l  = d\cdot\partial_ip_G(\mu) = d(d+1)p_{G}(\mu),
	\end{align*} where the last equality is by~(a).
\end{proof}

Our second lemma shows that $p_G(\mu)=\theta$ for every face-critical point $\mu$ of $p_G$ with $p_G(\mu)>0$ if we impose the following combinatorial condition on the cocycle $G$. We say that $G  \subseteq\binom{V}{d+1}$ is \emph{triangle-free} if for all $I \in \binom{V}{d}$ and distinct $i,j \in V \setminus I$ such that $I \cup \{i\}, I \cup \{j\} \in G$, the set $\{i,j\}$ is not a subset of any member of $G$. 

\begin{lemma}\label{lem:trianglefree}
Let  $G$ be a triangle-free $d$-cocycle of $\Delta(V)$, and let $\mu\in\Delta(V)$ with $p_G(\mu) > 0$ be face-critical for $p_G$. Then $p_G(\mu) = \theta$.
\end{lemma}

\begin{proof} We assume without loss of generality that $\supp \mu = V$. We prove the lemma by induction on the number of elements of $G$. In the base case, when  $G$ consists of a single $(d+1)$-element set, we have $|V|= d+1,$ and Lemma~\ref{lem:first}(a) implies that $\mu_i = \frac{1}{d+1}$ for every $i \in V$, and so $p_G(\mu) = \frac{(d+1)!}{(d+1)^{d+1}}.$
	
For the induction step, assume that $G$ contains at least two members. Then $|V| \geq d+2$, and as $G$ is a $d$-cocycle there exist	$I \in \binom{V}{d}$ and distinct $i,j \in V \setminus I$ such that $I \cup \{i\}, I \cup \{j\} \in G$.
As $G$ is triangle-free,   $\{i,j\}$ is not a subset of any member of $G$. 
 
For $I' \in \binom{V \setminus \{i,j\}}{d}$ we have $I' \cup \{i\}\in G$ if and only if $I' \cup \{j\} \in G$. Indeed, $I' \cup \{i,j\}$ contains an even number of sets in $G$ as $G$ is a cocycle, and the only subsets of this set that can belong to $G$ are $I' \cup \{i\}$ and $I' \cup \{j\}$.

Let $V'=V \setminus j$, let $G' \subseteq  \binom{V'}{d+1}$ be the restriction of $G$ to $V'$. As $G'$ is a restriction of $G$, $G'$ is a  triangle-free $d$-cocycle of $\Delta(V')$.   Let $\mu' \in \Delta(V')$ be defined by $\mu'_i = \mu_i + \mu_j$ and $\mu'_k =\mu_k$ for $k \in V' \setminus \{i\}$.  By the observation above we have $p_{G'}(\mu') = p_G(\mu)$ and $\partial_k p_{G'}(\mu') = \partial_k p_{G}(\mu)$ for every $k \in V'$. Thus $\mu'$ is face-critical for $p_{G'}$. By the induction hypothesis we have $p_G(\mu) = p_{G'}(\mu') = \theta$, as desired.
\end{proof}	

The following is the main result of this subsection.

\begin{lemma}\label{lem:index}
Let $G$ be a $d$-cocycle of $\Delta(V)$, and let $\mu\in\Delta(V)$ with $p_G(\mu) > 0$ be face-critical for $p_G$. Then the negative index of $p_G$ at $\mu$ is at least~$d$. If it equals $d$ then $p_G(\mu) = \theta$.
\end{lemma}

\begin{proof}
By replacing $V$ with $\supp(\mu)$	we assume without loss of generality that $V=\supp(\mu)$, i. e. $\mu_i>0$ for every $i\in V$. As $p_G(\mu)>0$ we have $G \neq \emptyset$ and so  $|V|\ge d+1$. Let $A=(\partial_i\partial_jp_G(\mu))_{i,j\in V}$ be the matrix of $\Hess p_G(\mu)$. Then $A$ is symmetric and has zero diagonal.

The proof of the lower bound on the negative index of $p_G$ has two main ingredients in addition to the preceding lemmas. Claim~\ref{cl:sos} writes the Hessian form plus $(d+1)p_G(\mu)$ times a positive definite form as a sum of squares, one for each member of $G$; after a change of variables, this bounds every eigenvalue on the tangent space below by $-(d+1)p_G(\mu)$. By Claim~\ref{cl:spectrum}, these eigenvalues sum to $-d(d+1)p_G(\mu)$, which forces at least $d$ of them to be negative. 

The following identity bounds the Hessian form from below.

\begin{claimnum}\label{cl:sos}
For every $x\in\R^V$,
\begin{equation*}
 x^{\mathsf T}Ax+(d+1)p_G(\mu)\sum_{i\in V}\frac{x_i^2}{\mu_i}=(d+1)!\sum_{e\in G}\Bigl(\prod_{l\in e}\mu_l\Bigr)\Bigl(\sum_{i\in e}\frac{x_i}{\mu_i}\Bigr)^2.
\end{equation*}
\end{claimnum}

\begin{proof}
	\begin{align*}
		(d+1)!&\sum_{e\in G}\Bigl(\prod_{l\in e}\mu_l\Bigr)\Bigl(\sum_{i\in e}\frac{x_i}{\mu_i}\Bigr)^2
	=(d+1)!\sum_{e\in G}\Bigl(\prod_{l\in e}\mu_l\Bigr)\Bigl(\sum_{i\in e}\frac{x_i^2}{\mu_i^2}+\sum_{\substack{i,j\in e\\ i\ne j}}\frac{x_ix_j}{\mu_i\mu_j}\Bigr)\\
		&=\sum_{i\in V}\frac{x_i^2}{\mu_i^2}\,(d+1)!\sum_{e\in G,\ e\ni i}\ \prod_{l\in e}\mu_l
	+\sum_{\substack{i,j\in V\\ i\ne j}}\frac{x_ix_j}{\mu_i\mu_j}\,(d+1)!\sum_{e\in G,\ e\supseteq\{i,j\}}\ \prod_{l\in e}\mu_l\\
		&=\sum_{i\in V}\frac{x_i^2}{\mu_i^2}\,\mu_i\,\partial_ip_G(\mu)+\sum_{\substack{i,j\in V\\ i\ne j}}\frac{x_ix_j}{\mu_i\mu_j}\,\mu_i\mu_jA_{ij}\\
		&=x^{\mathsf T}Ax+(d+1)p_G(\mu)\sum_{i\in V}\frac{x_i^2}{\mu_i}.
	\end{align*}
	The third equality follows from~\eqref{eq:derivbound} and~\eqref{eq:secondderiv}, multiplied by $\mu_i$ and by $\mu_i\mu_j$, respectively. The last equality uses Lemma~\ref{lem:first}(a) and the vanishing of the diagonal of~$A$.
\end{proof}

We next rescale the coordinates to eliminate the denominators $\mu_i$ on the left side. Let $\sqrt\mu = (\sqrt{\mu_i})_{i\in V} \in\R^V$ be a column vector and let $D =  \mathrm{diag} (\sqrt{\mu_i})_{i \in V}$ be the corresponding diagonal matrix. Let $\widetilde A = DAD$, so that
\begin{equation*}
 \widetilde A_{ij}=\sqrt{\mu_i\mu_j}\,A_{ij},\qquad \sqrt\mu=(\sqrt{\mu_i})_{i\in V}.
\end{equation*}
For $z\in\R^V$, let $x\in\R^V$ have coordinates $x_i=\sqrt{\mu_i}\,z_i$. Then
\begin{equation}\label{eq:subst}
 x^{\mathsf T}Ax=z^{\mathsf T}\widetilde Az,\qquad \sum_{i\in V}\frac{x_i^2}{\mu_i}=\|z\|^2,\qquad \sum_{i\in V}x_i=\sqrt\mu\cdot z,\qquad \sum_{i\in e}\frac{x_i}{\mu_i}=\sum_{i\in e}\frac{z_i}{\sqrt{\mu_i}},
\end{equation}
where $\cdot$ and $\|\cdot\|$ denote the standard inner product and norm. In particular, $z\mapsto x$ is a linear bijection from $(\sqrt\mu)^\perp=\{z\in\R^V:\sqrt\mu\cdot z=0\}$ onto $\R^V_0$ that carries the form $z^{\mathsf T}\widetilde Az$ to the Hessian form, and Claim~\ref{cl:sos} becomes
\begin{equation}\label{eq:sosz}
 z^{\mathsf T}\widetilde Az+(d+1)p_G(\mu)\,\|z\|^2=(d+1)!\sum_{e\in G}\Bigl(\prod_{l\in e}\mu_l\Bigr)\Bigl(\sum_{i\in e}\frac{z_i}{\sqrt{\mu_i}}\Bigr)^2\qquad(z\in\R^V).
\end{equation}
Moreover, $\|\sqrt\mu\|^2=\sum_i\mu_i=1$; $\tr\widetilde A=0$, since $A$ has zero diagonal; and by Lemma~\ref{lem:first}(b), $(\widetilde A\sqrt\mu)_i=\sqrt{\mu_i}\sum_jA_{ij}\mu_j=d(d+1)p_G(\mu)\sqrt{\mu_i}$ for every $i$, that is,
\begin{equation}\label{eq:eigen}
 \widetilde A\sqrt\mu=d(d+1)p_G(\mu)\,\sqrt\mu.
\end{equation}
Hence $\widetilde A$ maps $(\sqrt\mu)^\perp$ into itself: if $\sqrt\mu\cdot z=0$, then $\sqrt\mu\cdot\widetilde Az=(\widetilde A\sqrt\mu)\cdot z=0$. The restriction of $\widetilde A$ to $(\sqrt\mu)^\perp$ is self-adjoint, since $\widetilde A$ is symmetric. By the real spectral theorem~\cite[\S7B]{Axler}, the space $(\sqrt\mu)^\perp$, which has dimension $|V|-1$, has an orthonormal basis $y_1,\ldots,y_{|V|-1}$ consisting of eigenvectors of  $\widetilde A$ with corresponding eigenvalues $\lambda_1\le\cdots\le\lambda_{|V|-1}$. We write $y_k(i)$ for the $i$th coordinate of $y_k$.

\begin{claimnum}\label{cl:count}
The negative index of $p_G$ at $\mu$ is the number of $k$ with $\lambda_k<0$.
\end{claimnum}

\begin{proof}
Let $s$ be this number, so that $\lambda_k<0$ exactly for $k\le s$. By~\eqref{eq:subst}, the negative index is the largest dimension of a subspace $L\subseteq(\sqrt\mu)^\perp$ such that $z^{\mathsf T}\widetilde Az<0$ for every nonzero $z\in L$. For $z=\sum_kc_ky_k$, orthonormality gives $z^{\mathsf T}\widetilde Az=\sum_k\lambda_kc_k^2$. This is negative for nonzero $z$ in the span of $y_1,\ldots,y_s$, so the negative index is at least $s$. Conversely, let $L\subseteq(\sqrt\mu)^\perp$ be a subspace with $\dim L>s$. The dimensions of $L$ and of the span of $y_{s+1},\ldots,y_{|V|-1}$ add up to more than $|V|-1$, so these subspaces share a nonzero vector $z=\sum_{k>s}c_ky_k$, and $z^{\mathsf T}\widetilde Az=\sum_{k>s}\lambda_kc_k^2\ge0$.
\end{proof}

The next claim establishes the desired lower bound on the negative index of $p_G$ by bounding the eigenvalues of the restriction of $\widetilde A$ to  $(\sqrt\mu)^\perp$.

\begin{claimnum}\label{cl:spectrum}
\begin{enumerate}[label=\textup{(\alph*)}]
\item For every $k$, $\lambda_k\ge-(d+1)p_G(\mu)$, with equality if and only if $\sum_{i\in e}y_k(i)/\sqrt{\mu_i}=0$ for every $e\in G$.
\item $\sum_k\lambda_k=-d(d+1)p_G(\mu)$.
\end{enumerate}
\end{claimnum}

\begin{proof}
(a) Since $\|y_k\|=1$, we have $\lambda_k=y_k^{\mathsf T}\widetilde Ay_k$, and~\eqref{eq:sosz} with $z=y_k$ gives
\begin{equation*}
 \lambda_k+(d+1)p_G(\mu)=(d+1)!\sum_{e\in G}\Bigl(\prod_{l\in e}\mu_l\Bigr)\Bigl(\sum_{i\in e}\frac{y_k(i)}{\sqrt{\mu_i}}\Bigr)^2.
\end{equation*}
The right-hand side is a sum of squares with positive coefficients; it is nonnegative, and it vanishes exactly when every inner sum vanishes.

(b) The vectors $\sqrt\mu,y_1,\ldots,y_{|V|-1}$ form an orthonormal basis of $\R^V$; let $F$ be the orthogonal matrix with these columns. By~\eqref{eq:eigen} and the choice of $y_1,\ldots,y_{|V|-1}$, the matrix $F^{\mathsf T}\widetilde AF$ is diagonal with diagonal entries $d(d+1)p_G(\mu),\lambda_1,\ldots,\lambda_{|V|-1}$. Its trace is $\tr(\widetilde AFF^{\mathsf T})=\tr\widetilde A=0$.
\end{proof}

Claim~\ref{cl:spectrum} immediately implies that at least $d$ of the numbers $\lambda_k$ are negative. Thus by Claim~\ref{cl:count} the negative index of $p_G$ is at least $d$. This establishes the first claim of the lemma.

Assume for the rest of the proof that the negative index equals $d$, i.e. exactly $d$ eigenvalues $\lambda_k$ are negative. By Claim~\ref{cl:spectrum} we have the following.

\begin{claimnum}\label{cl:equal}
	We have $\lambda_k=-(d+1)p_G(\mu)$ for $1\le k\le d$ and $\lambda_k=0$ for $k>d$. Moreover, for every $f\in G$ and $1\le k\le d$,
	\begin{equation*}
		\sum_{v\in f}\frac{y_k(v)}{\sqrt{\mu_v}}=0.
	\end{equation*}
\end{claimnum}

By Lemma~\ref{lem:trianglefree} the following last claim implies $p_G(\mu)=\theta$ and thus the current lemma.

\begin{claimnum}\label{cl:exchangezero}
	$G$ is triangle-free.
\end{claimnum}

\begin{proof}
	Let $I \in \binom{V}{d}$ and distinct $i,j \in V \setminus I$ such that $I \cup \{i\}, I \cup \{j\} \in G$. Our goal is to show that $G$ has no members containing both $i$ and $j$. Subtracting the identities in Claim~\ref{cl:equal} for  $I \cup \{i\}$ and $I \cup \{j\}$ we get
	\begin{equation*}
		\frac{y_k(i)}{\sqrt{\mu_i}}
		=\frac{y_k(j)}{\sqrt{\mu_j}}
		\qquad(1\le k\le d).
	\end{equation*}
	Let $e_i,e_j\in\R^V$ be the corresponding standard basis vectors, and let $z=D^{-1}(e_i-e_j)$. Then
	\begin{equation*}
		\sqrt\mu\cdot z=1-1=0,\qquad
		y_k\cdot z=
		\frac{y_k(i)}{\sqrt{\mu_i}}
		-\frac{y_k(j)}{\sqrt{\mu_j}}=0
		\qquad(1\le k\le d).
	\end{equation*}
	Since $\sqrt\mu,y_1,\ldots,y_{|V|-1}$ form an orthonormal basis, $z$ lies in the span of the vectors $y_k$ with $k>d$. Their eigenvalues are zero by Claim~\ref{cl:equal}, so $\widetilde Az=0$. As $\widetilde A=DAD$ and $A$ is symmetric with zero diagonal,
	\begin{equation*}
		0=z^{\mathsf T}\widetilde Az
		=(e_i-e_j)^{\mathsf T}A(e_i-e_j)
		=-2A_{ij}.
	\end{equation*}
	By~\eqref{eq:secondderiv}, $A_{ij}$ is a sum of positive summands over the members of $G$ containing $i$ and $j$. As $A_{ij}=0$, there are no such members.
\end{proof}
\end{proof}

\subsection{Morse theory}\label{sec:morsetheory}

In this subsection we derive Theorem~\ref{thm:weights} from Lemma~\ref{lem:index} using Morse theory. 

Morse theory provides a way of deriving topological properties of a manifold by studying smooth functions on it. Milnor~\cite{Milnor} serves as our primary reference. We only need smooth manifolds $M\subseteq\R^k$, and our application takes $M$ to be a sphere. For smooth maps, diffeomorphisms, tangent spaces $T_xM\subseteq\R^k$, and differentials, we refer to~\cite[\S1]{MilnorTDV}. A \emph{chart} of $M$ is the inverse of a parametrization of an open subset of $M$, called a system of coordinates in~\cite[\S1]{MilnorTDV}. For a smooth function $h:M\to\R$, we write $dh_x:T_xM\to\R$ for its differential at $x\in M$; a \emph{critical point} of $h$ is a point $x$ with $dh_x=0$, and a \emph{critical value} is the value of $h$ at a critical point.

Let $M\subseteq\R^k$ be a compact smooth manifold without boundary and $f:M\to\R$ a smooth function. At a critical point of $f$, the Hessian is the quadratic form on the tangent space given by the second derivatives in a chart; its \emph{negative index}, the largest dimension of a subspace on which it is negative definite, is independent of the chart~\cite[\S2]{Milnor}. The function $f$ is \emph{Morse} if all these Hessians are nondegenerate; its critical points are then isolated~\cite[Corollary~2.3]{Milnor}, hence finitely many. A value is \emph{regular} if it is not a critical value. We fix finitely many charts of $M$ and compact subsets of their domains that cover $M$, and say that a smooth function $g$ is \emph{$\varepsilon$-close} to $f$ if $g-f$ and its partial derivatives of order at most two in these charts have absolute value less than $\varepsilon$ on these compact sets. We use the following standard Morse theory results, in the generality we need.

\begin{theorem}\label{thm:morse}
Let $f:M\to\R$ be smooth.
\begin{enumerate}[label=\textup{(\alph*)}]
\item If $a<b$ and $f^{-1}([a,b])$ contains no critical point, then $\{f\le a\}$ is a strong deformation retract of $\{f\le b\}$~\cite[Theorem~3.1 and its proof, pp.~12--13]{Milnor}.
\item If $f$ is Morse and $p$ is the only critical point with value $v=f(p)$, then $H_i(\{f\le v+\varepsilon\},\{f\le v-\varepsilon\})=0$ for every sufficiently small $\varepsilon>0$ and every $i$ different from the negative index of $f$ at~$p$~\cite[Theorem~5.3.1, p.~278]{Goresky}.
\item For every $\varepsilon>0$ there is a Morse function that is $\varepsilon$-close to $f$~\cite[Corollary~6.8, p.~37]{Milnor}.
\end{enumerate}
\end{theorem}

Theorem~\ref{thm:morse}(a) and (b) appear to provide a perfect toolkit for proving Theorem~\ref{thm:weights} with Lemma~\ref{lem:index} used to control the index of critical points. There are however a couple of technical issues: \begin{itemize} \item The simplex $\Delta(V)$ is not a smooth manifold without boundary,
	\item the function $p_G$ need not be Morse in the band we are interested in, as we have no control over non-degeneracy of critical points with values in $(0,\theta)$.	\end{itemize}
The rest of this section is devoted to overcoming these issues, by pulling back $p_G$ from $\Delta(V)$ to a sphere and using Theorem~\ref{thm:morse}(c) to apply Morse theory to non-Morse smooth functions.

\begin{lemma}\label{lem:band}
Let $\varphi:M\to\R$ be smooth, let $a<b$ and $\eta>0$, and suppose that every critical point $x$ of $\varphi$ with $a\le\varphi(x)\le b+\eta$ has negative index at least $d+1$. Let $0\le j\le d$ and $c\in S_j(\{\varphi\le b\})$ be a singular chain with $\partial c\in S_{j-1}(\{\varphi\le a\})$. Then $c=\partial w+z$ for some $w\in S_{j+1}(\{\varphi<b+\eta\})$ and $z\in S_j(\{\varphi<a+\eta\})$.
\end{lemma}

\begin{proof}
We first show that there is $\varepsilon\in(0,\eta/3]$ such that, for every function $h$ that is $\varepsilon$-close to $\varphi$, every critical point $x$ of $h$ with $a\le h(x)\le b+\eta$ has negative index at least $d+1$. Indeed, otherwise there are functions $h_k$ that are $(1/k)$-close to $\varphi$ and critical points $x_k$ of $h_k$ with $a\le h_k(x_k)\le b+\eta$ and negative index at most $d$. A subsequence of $(x_k)$ lies in one of the fixed compact sets and converges to a point $x$, and in the corresponding chart the derivatives of order at most two of $h_k$ at $x_k$ converge to those of $\varphi$ at $x$. Thus $x$ is a critical point of $\varphi$ with $a\le\varphi(x)\le b+\eta$, so its Hessian matrix is negative definite on a $(d+1)$-dimensional subspace, and so is the Hessian matrix of $h_k$ at $x_k$ for large $k$, a contradiction.

By Theorem~\ref{thm:morse}(c), let $h$ be a Morse function that is $\varepsilon$-close to $\varphi$, and choose regular values $a'\in(a+\varepsilon,a+2\varepsilon)$ and $b'\in(b+\varepsilon,b+2\varepsilon)$ of $h$ with $a'<b'$. Since $|h-\varphi|<\varepsilon$ and $3\varepsilon\le\eta$,
\begin{equation*}
 \{\varphi\le a\}\subseteq\{h\le a'\}\subseteq\{\varphi<a+\eta\},\qquad
 \{\varphi\le b\}\subseteq\{h\le b'\}\subseteq\{\varphi<b+\eta\}.
\end{equation*}
By the characterization of relative homology in ``\nameref{sec:chains}'', it suffices to show $H_j(\{h\le b'\},\{h\le a'\})=0$. By the first paragraph, every critical point of $h$ with value in $[a',b']$ has negative index at least $d+1$.

Let $p_1,\ldots,p_r$ be these critical points; their values lie in $(a',b')$. We first make their values distinct. Choose pairwise disjoint open sets $U_i\ni p_i$ that contain no other critical point and have closures in $h^{-1}((a',b'))$, and smooth functions $\omega_i:M\to[0,1]$ that vanish outside a compact subset of $U_i$ and are equal to $1$ near~$p_i$. Let $K$ be the closure of the set of points at which the differential of some $\omega_i$ is nonzero. Then $K$ is compact and contains no critical point of $h$, as it is contained in $\bigcup_iU_i$ and avoids a neighborhood of every~$p_i$. Thus the norm of $dh_x$ with respect to the Euclidean inner product on $T_xM$ is bounded away from zero for $x\in K$. Hence, for constants $u_i$ of sufficiently small absolute value, $h'=h+\sum_iu_i\omega_i$ has the same critical points as $h$, with the same Hessians, and takes values in $(a',b')$ on every $U_i$, so that $\{h'\le a'\}=\{h\le a'\}$ and $\{h'\le b'\}=\{h\le b'\}$. We choose the $u_i$ so that the values $h'(p_i)=h(p_i)+u_i$ are distinct, and order them as $v_1<\cdots<v_r$.

If $r=0$, Theorem~\ref{thm:morse}(a) and Lemma~\ref{lem:prism}(b) give the claim. Otherwise choose $a'=t_0<v_1<t_1<\cdots<v_r<t_r=b'$. By Theorem~\ref{thm:morse}(b), there are $\varepsilon_i>0$ with $t_{i-1}<v_i-\varepsilon_i<v_i+\varepsilon_i<t_i$ and $H_j(\{h'\le v_i+\varepsilon_i\},\{h'\le v_i-\varepsilon_i\})=0$, since $j\le d$ is smaller than the negative index at the critical point with value~$v_i$. By Theorem~\ref{thm:morse}(a) and Lemma~\ref{lem:prism}(b), the relative homology of the pairs $(\{h'\le v_i-\varepsilon_i\},\{h'\le t_{i-1}\})$ and $(\{h'\le t_i\},\{h'\le v_i+\varepsilon_i\})$ vanishes in all degrees. Lemma~\ref{lem:prism}(c) now gives $H_j(\{h'\le b'\},\{h'\le a'\})=0$.
\end{proof}

We are now ready to prove Theorem~\ref{thm:weights}.

\begin{proof}
Put $B=\{p_G<\theta\}$ and $I=\{p_G=0\}=|\Ind(G)|$. If $G=\varnothing$, then $B=I$, and there is nothing to prove. Otherwise put $N=|V|\ge d+1$, let $S^{N-1}\subset\R^V$ be the unit sphere, and define
\begin{equation*}
 q:S^{N-1}\to\Delta(V),\quad q(x)=(x_i^2)_{i\in V},\qquad F=p_G\circ q,\qquad
 s:\Delta(V)\to S^{N-1},\quad s(\mu)=(\sqrt{\mu_i})_{i\in V}.
\end{equation*}
Then $F$ is smooth, $s$ is continuous, and $q\circ s$ is the identity, so $F\circ s=p_G$.

Let $x$ be a critical point of $F$ with $0<F(x)<\theta$, and put $\mu=q(x)$ and $\sigma=\supp\mu$. The points $y$ of the coordinate sphere $S^{N-1}\cap\R^\sigma$ with $y_ix_i>0$ for all $i\in\sigma$ form an open subset $P$ of it containing $x$, and $q$ maps $P$ diffeomorphically onto $\Delta^\circ(\sigma)$. Hence $\mu$ is a critical point of the restriction of $p_G$ to $\Delta^\circ(\sigma)$, that is, $\mu$ is face-critical for $p_G$. Since first derivatives vanish at critical points, the chain rule shows that the Hessian of $F$ at $x$, restricted to $T_xP$, is carried by the differential of $q$ to the Hessian of $p_G$ at $\mu$, restricted to $\R^\sigma_0$. As $p_G(\mu)<\theta$, by Lemma~\ref{lem:index} the  Hessian of $p_G$ at $\mu$  is negative definite on a subspace of dimension $d+1$. Thus the negative index of $F$ at $x$ is at least $d+1$.

Let $0\le j\le d$ and $c\in S_j(B)$ with $\partial c\in S_{j-1}(I)$. Let $\alpha=(d+1)!/(2N)^{d+1}$; since $N\ge d+1$, we have $2^{d+1}\alpha\le\theta$. Choose $b\in(\alpha/3,\theta)$ with $c \in S_j(\{p_G\le b\})$, and put $\eta=\min\{\alpha/3,(\theta-b)/2\}$. The chain $s_\#c$ lies in $S_j\{F\le b\}$, and its boundary lies in $s(I)\subseteq\{F=0\}$. Lemma~\ref{lem:band} applies to $F$ with $a=\alpha/3$, and applying $q_\#$ to the resulting identity gives
\begin{equation*}
 c=\partial w+y,\qquad w\in S_{j+1}(\{p_G<b+\eta\})\subseteq S_{j+1}(B),\quad y\in S_j(\{p_G<\alpha\}),
\end{equation*}
with $\partial y=\partial c\in S_{j-1}(I)$.

It remains to deform $\{p_G<\alpha\}$ into $I$ inside $B$, keeping $I$ inside $I$. For $\mu\in\Delta(V)$ and $t\in[0,1]$ let $\Theta(\mu,t)$ be the vector with coordinates $\max\{\mu_i-t/(2N),0\}$, divided by the sum of these coordinates, which is at least $1-t/2\ge1/2$. Then $\Theta$ is continuous, $\Theta(\mu,0)=\mu$, and $\supp\Theta(\mu,t)\subseteq\supp\mu$, so $\Theta$ maps $I\times[0,1]$ into $I$, as $\Ind(G)$ is closed under taking nonempty subsets. Since $p_G$ is homogeneous of degree $d+1$ and nondecreasing in each coordinate on $[0,\infty)^V$, we have $p_G(\Theta(\mu,t))\le2^{d+1}p_G(\mu)$, so $\Theta$ maps $\{p_G<\alpha\}\times[0,1]$ into $B$. If $p_G(\mu)<\alpha$, then every $e\in G$ contains some $i$ with $\mu_i<1/(2N)$, because $(d+1)!\prod_{i\in e}\mu_i\le p_G(\mu)<(d+1)!\,(2N)^{-(d+1)}$. Thus $\Theta(\mu,1)_i=0$, so $\Theta(\mu,1)\in I$. Lemma~\ref{lem:prism}(a), applied to $\Theta$ on $\{p_G<\alpha\}\times[0,1]$, gives
\begin{equation*}
 y=\partial\Pi(y)+\Pi(\partial y)+\Theta(\cdot,1)_\#y
\end{equation*}
with $\Pi(y)\in S_{j+1}(B)$ and $\Pi(\partial y)+\Theta(\cdot,1)_\#y\in S_j(I)$. Hence $c=\partial(w+\Pi(y))+\Pi(\partial y)+\Theta(\cdot,1)_\#y$ is the required decomposition.
\end{proof}

\section{Concluding remarks}\label{sec:remarks}

\subsubsection*{Maps to manifolds.}

The proof of Lemma~\ref{lem:intersection} following the argument in~\cite{DKW} extends to continuous maps $f : \Delta^n \to M$ for any PL $d$-manifold $M$. As Lemma~\ref{lem:intersection} is the only place where the map enters the proof, Theorems~\ref{thm:finite} and~\ref{thm:main} extend to maps to PL $d$-manifolds.

\subsubsection*{Brambles.}

Following~\cite{Norin} let a \emph{pure $d$-dimensional bramble} on $\Delta^n$ be a collection $\mc{B}$ of pure $d$-dimensional subcomplexes of $\Delta^n$ such that \begin{itemize} \item $\mc{B}$ is closed under unions, and \item $\tilde{H}^{i}(Z) = 0$ for every $Z \in \mc{B}$ and $0 \leq i \leq d-1$.\end{itemize}
$d$-dimensional brambles are high-dimensional analogues of graph brambles, which Seymour and Thomas~\cite{ST} introduced under the name \emph{screens}.

\emph{The order $\brm{ord}(\mc{B})$ of the bramble $\mc{B}$} is the minimum size of a collection $S \subseteq \binom{[n]}{d+1}$ of $d$-dimensional faces of $\Delta^n$ such that every element of $\mc{B}$ contains some face in $S$. 
It follows from the results of~\cite{Norin} that \begin{equation}\label{e:bramble} \ov(f) \geq \brm{ord}(\mc{B})\end{equation}
for every continuous $f: \Delta^n \to \R^d$ and every pure $d$-dimensional bramble $\mc{B}$ on $\Delta^n$.

The bound \eqref{e:bramble} initially motivated the investigation that led to the results of this paper, but our resulting overlap bounds do not explicitly use it. 

It seems much harder and likely impossible to obtain sharp topological overlap bounds using only \eqref{e:bramble}, although it can be done for $d=2$. It would be interesting to determine the limits of the method.

\begin{question}
What are the asymptotics of the maximum order of a pure $d$-dimensional bramble on $\Delta^n$? 	
\end{question}

\appendix

\section{A running example}\label{app:example}

In this appendix we illustrate the proof of Theorem~\ref{thm:finite} on a small example: Appendix~\ref{app:example-cochains} describes the output of Lemma~\ref{lem:intersection} for an affine map $f:\Delta^4\to\R^2$, and Appendix~\ref{app:example-regions} describes the sets of Section~\ref{sec:proof} for this map.

\subsection{Intersection cochains}\label{app:example-cochains}

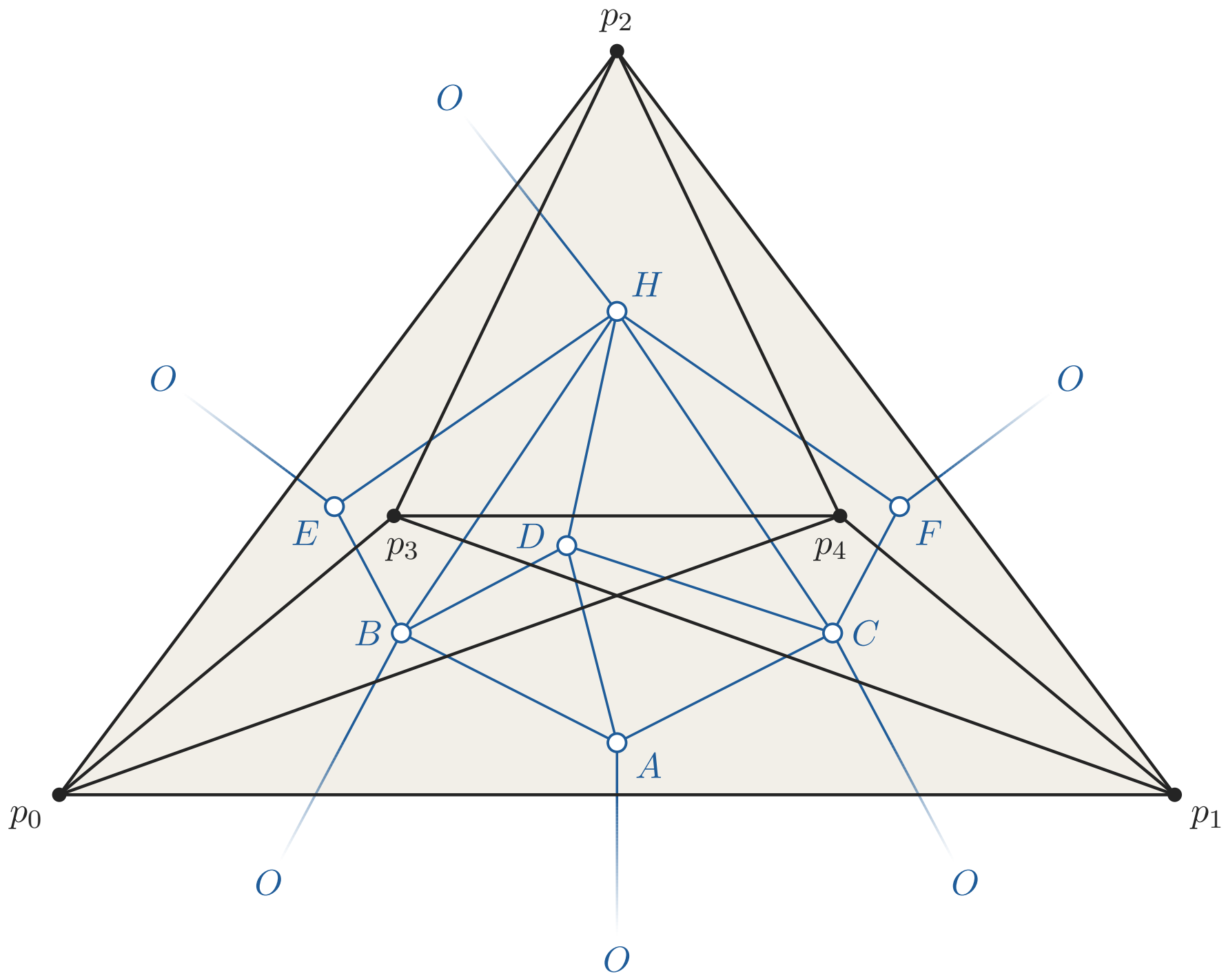
\begin{figure}[H]
\centering
\resizebox{0.75\textwidth}{!}{% Extracted from Experiments/five_point_T.tex.
% Requires \usepackage{tikz} and \usetikzlibrary{calc} in the host preamble.
% The target sphere T of the five-point example, drawn on the compactified
% plane R^2 \cup \{\infty\} with the vertex O at infinity, together with the
% image f(|\Delta_5^{(2)}|) (Experiments/five_point_example.md).
% Edges of T ending at O are drawn as rays fading towards infinity.
%
% For legibility the interior points are drawn at p_4=(-6/5,3/2), p_5=(6/5,3/2)
% instead of (-1,1), (1,1); then X = p_1p_5 \cap p_2p_4 = (0,15/14).  The
% drawing is combinatorially identical: same unique crossing X, same seven
% bounded regions, and every crossing set a(e), every vertex cochain a(v) and
% every top cochain a(tau) agrees with the tables in the note.
\begingroup
\definecolor{srcfill}{RGB}{242,239,232}
\definecolor{srcline}{RGB}{36,36,36}
\definecolor{tcol}{RGB}{31,92,153}

\begin{tikzpicture}[
  scale=2,
  sedge/.style={srcline, line width=0.95pt, line join=round, line cap=round},
  tedge/.style={tcol, line width=0.75pt, line cap=round, line join=round},
  svert/.style={circle, fill=srcline, inner sep=0pt, minimum size=4.4pt},
  tvert/.style={circle, fill=white, draw=tcol, line width=0.8pt, inner sep=0pt, minimum size=5.6pt},
  slabel/.style={text=srcline, inner sep=1.2pt},
  tlabel/.style={text=tcol, inner sep=1.2pt},
  % edge from a finite vertex of T to O = infinity, solid until it has left
  % the source triangle, then fading:
  %   \pic{toinfty={<vertex>}{<angle>}{<solid length>}{<fade length>}};
  pics/toinfty/.style n args={4}{code={
    \def\S{#3}\def\F{#4}\def\N{40}
    \draw[tedge] (#1) -- ($(#1)+(#2:\S)$);
    \foreach \k in {1,...,\N}{
      \pgfmathsetmacro{\a}{\S+\F*(\k-1)/\N}
      \pgfmathsetmacro{\b}{\S+\F*\k/\N}
      \pgfmathsetmacro{\op}{(1-(\k-0.5)/\N)^1.25}
      \draw[tedge, line cap=butt, draw opacity=\op] ($(#1)+(#2:\a)$) -- ($(#1)+(#2:\b)$);
    }
    \node[tlabel] at ($(#1)+(#2:{\S+\F+0.25})$) {$O$};
  }},
]
  % ---------------- image of Delta_5^{(2)} ----------------
  \coordinate (p1) at (-3,0);
  \coordinate (p2) at (3,0);
  \coordinate (p3) at (0,4);
  \coordinate (p4) at (-1.2,1.5);
  \coordinate (p5) at (1.2,1.5);
  \fill[srcfill] (p1) -- (p2) -- (p3) -- cycle;

  % ---------------- target sphere T (O at infinity) ----------------
  \coordinate (A) at (0,0.28);      % in 12X
  \coordinate (B) at (-1.16,0.87);  % in 14X
  \coordinate (C) at (1.16,0.87);   % in 25X
  \coordinate (D) at (-0.27,1.34);  % in 45X
  \coordinate (E) at (-1.52,1.55);  % in 134
  \coordinate (F) at (1.52,1.55);   % in 235
  \coordinate (H) at (0,2.6);       % in 345

  % edges to O
  \pic {toinfty={A}{270}{0.58}{1.5}};     % crosses 12
  \pic {toinfty={B}{242}{1.29}{1.5}};     % crosses 15, 12
  \pic {toinfty={C}{298}{1.29}{1.5}};     % crosses 24, 12
  \pic {toinfty={E}{143.13}{0.55}{1.5}};  % crosses 13
  \pic {toinfty={F}{36.87}{0.55}{1.5}};   % crosses 23
  \pic {toinfty={H}{128}{1.17}{1.5}};     % crosses 34, 13
  % remaining cube edges
  \draw[tedge] (A) -- (B)  (A) -- (C)  (B) -- (D)  (C) -- (D)
               (B) -- (E)  (C) -- (F)  (E) -- (H)  (F) -- (H)  (D) -- (H);
  % remaining diagonals (AD passes to the left of X)
  \draw[tedge] (B) -- (H)  (C) -- (H)  (A) -- (D);

  % source edges drawn over the target edges
  \draw[sedge] (p1) -- (p2) -- (p3) -- cycle;
  \draw[sedge] (p1) -- (p4) -- (p3)  (p2) -- (p5) -- (p3)  (p4) -- (p5);
  \draw[sedge] (p1) -- (p5)  (p2) -- (p4);

  \foreach \p in {p1,p2,p3,p4,p5} \node[svert] at (\p) {};
  \foreach \v in {A,B,C,D,E,F,H} \node[tvert] at (\v) {};

  % ---------------- labels ----------------
  \node[slabel] at ($(p1)+(215:0.216)$) {$p_0$};
  \node[slabel] at ($(p2)+(325:0.216)$) {$p_1$};
  \node[slabel] at ($(p3)+(90:0.163)$)  {$p_2$};
  \node[slabel] at ($(p4)+(285:0.188)$) {$p_3$};
  \node[slabel] at ($(p5)+(255:0.188)$) {$p_4$};
  \node[tlabel] at ($(A)+(324:0.21)$)   {$A$};
  \node[tlabel] at ($(B)+(180:0.178)$)  {$B$};
  \node[tlabel] at ($(C)+(0:0.178)$)    {$C$};
  \node[tlabel] at ($(D)+(165:0.202)$)  {$D$};
  \node[tlabel] at ($(E)+(222:0.21)$)   {$E$};
  \node[tlabel] at ($(F)+(318:0.21)$)   {$F$};
  \node[tlabel] at ($(H)+(42:0.218)$)   {$H$};
\end{tikzpicture}
\endgroup}
\caption{An example of an affine map $f: \Delta^4 \to \R^2$ and a simplicial complex $T$ satisfying the conditions of Lemma~\ref{lem:intersection}.}
\label{fig:five-point-T}
\end{figure}

We use the following example, pictured in Figure~\ref{fig:five-point-T}, to present some of the ideas of the proof in a concrete setting. Let $f: \Delta^4 \to \R^2$ be an affine map with the vertices of $\Delta^4$
mapped to $\{p_0,\ldots,p_4\}$. The simplicial complex $T$ is a triangulation of the $2$-skeleton of the $3$-dimensional cube represented in blue with the vertex $O$ at infinity. Apply the intersection-cochain construction from the proof of Lemma~\ref{lem:intersection} directly to $f$. The resulting map $a:C_*(T)\to C^{d-*}(\Delta^4)$ can be read off from the figure. For a triangle $\tau$ of $T$, the support of $a(\tau)$ consists of those $i$ for which $p_i$ lies in the region represented by $\tau$. Each such region contains at most one $p_i$. For example, $\supp a(OCF) = \{1\}$ and $a(ACD)  = 0$. For an edge $\sigma$ of $T$, $a(\sigma)$ is supported exactly on the edges of $\Delta^4$ the images of which cross $\sigma$, e.g. $\supp a(AD)= \{04,13\}$.
Finally, for a vertex $v$ of $T$, $a(v)$ is supported on triangles of $\Delta^4$ the  images of which  contain $v$, e.g. $\supp a(D) =  \{012, 024, 034, 123, 134\}$. 

It is not hard to see that the map $a$ defined in this way is a chain--cochain map. By linearity, it suffices to check that $\delta a(\sigma)=a(\partial \sigma)$ for every simplex $\sigma$ of $T$, e.g. when $\sigma =\{u,v\}$ is an edge this amounts to checking that the symmetric difference of  $\supp a(u)$ and $\supp a(v)$ consists exactly of the triangles of $\Delta^{4}$ which contain an odd number of edges crossing $uv$. 

The sum $\zeta$ of all triangles of $T$ is a $2$-cycle of $T$ and satisfies $a(\zeta)=\one$.

\subsection[The sets S, I, and B]{The sets $S_\tau$, $I_\tau$, and $B_\tau$}\label{app:example-regions}

We need to slightly modify the construction  in  Section~\ref{sec:proof} to make it non-trivial on our small example.
With the sets $S_\tau$ defined before Lemma~\ref{lem:covers}, $|V_\tau|\le2$, and hence $G_\tau=\varnothing$, for every triangle $\tau$ of $T$. For each simplex $\eta$ of $T$ of positive dimension, instead of the vertex set of a maximal family of pairwise vertex-disjoint members of $\supp a(\eta)$ we take a smallest set meeting every member of $\supp a(\eta)$, and define $S_\tau$ to be the union of these sets over the faces of $\tau$ of positive dimension. For a triangle $\eta$ nothing changes, as $\supp a(\eta)$ has at most one element by Lemma~\ref{lem:intersection}\ref{I:gp}. For an edge $e\ne AD$, some vertex lies in every member of $\supp a(e)$, and we take one such vertex, choosing $4$, or else $3$, when possible; for $AD$ we take $\{3,4\}$. For instance, $S_{BD}=S_{BHE}=\{3\}$, $S_{CD}=S_{CDH}=S_{CHF}=\{4\}$, and $S_{AD}=S_{ABD}=S_{BDH}=\{3,4\}$.

These sets retain the nesting, vanishing, and relative-homology properties used in the proof. However, the uniform vector does not belong to every $B_\tau$: indeed,
\[
 p_{G_D}(b_{[4]})=\frac{6}{25}>\frac{2}{9}=\theta,
\]
so $b_{[4]}\notin B_D$.

Since $a(O)=0$, we have $G_\tau=\varnothing$ and $I_\tau=\Delta(V_\tau)$ whenever $O\in\tau$. Every other $G_\tau$ except $G_D$ is complete tripartite: $V_\tau$ is partitioned into nonempty parts $P_1,P_2,P_3$, and $G_\tau$ consists of the $3$-element sets meeting all three parts. Then a face of $\Delta(V_\tau)$ belongs to $I_\tau$ if and only if it misses a part, so $I_\tau=\Delta(P_1\cup P_2)\cup\Delta(P_1\cup P_3)\cup\Delta(P_2\cup P_3)$. For instance, $G_{BD}=G_{BHE}=\{012,024\}$ on $0124$ has parts $0,2,14$, $G_{CD}=G_{CDH}=G_{CHF}=\{012,123\}$ on $0123$ has parts $1,2,03$, and $G_{AD}=G_{ABD}=G_{BDH}=\{012\}$ on $012$; thus
\begin{equation*}
 I_{BD}=\Delta(014)\cup\Delta(124)\cup\Delta(02),\qquad
 I_{CD}=\Delta(013)\cup\Delta(023)\cup\Delta(12),\qquad
 I_{AD}=\partial\Delta(012).
\end{equation*}
Each of $I_{BD}$ and $I_{CD}$ consists of a disk formed by two triangles, together with an additional edge joining opposite corners, so it has the homotopy type of a circle. They lie in different tetrahedra and meet in $I_{AD}$.

The exception is $G_D=\supp a(D)$, which has five members. The complex $I_D$ consists of the triangles $013,014,124,234,023$ and their faces. Consecutive triangles in this cyclic order share the edges $01,14,24,23,03$, and the other five edges form the single boundary cycle $0\,2\,1\,3\,4\,0$, so $|I_D|$ is a M\"obius band. 
 
Let $$B'_{\tau} =  B_\tau\cap\Delta(V_\tau)=\{\mu\in\Delta(V_\tau):p_{G_\tau}(\mu)<\theta\}.$$ The proof of Lemma~\ref{lem:carriers}(iii) shows that $B'_{\tau}$ is a strong deformation retract of~$B_\tau$. If $G_\tau$ is complete tripartite, then $p_{G_\tau}(\mu)=6\,\mu(P_1)\mu(P_2)\mu(P_3)\le6/27=\theta$ for $\mu\in\Delta(V_\tau)$ with equality if and only if $\mu(P_1)=\mu(P_2)=\mu(P_3)=1/3$. Such $\mu$ form $\conv\{b_\sigma:\sigma\in G_\tau\}$, since every $b_\sigma$ satisfies these equations and every such $\mu$ is the convex combination $\sum27\,\mu_i\mu_j\mu_k\,b_{\{i,j,k\}}$ over $i\in P_1$, $j\in P_2$, $k\in P_3$. Hence
\begin{equation*}
 B'_\tau=\Delta(V_\tau)\setminus\conv\{b_\sigma:\sigma\in G_\tau\}.
\end{equation*}
Thus $B'_{BD}$ is $\Delta(0124)$ with a segment from $b_{012}$ to $b_{024}$ removed, $B'_{CD}$ is $\Delta(0123)$ with the segment from $b_{012}$ to $b_{123}$ removed, and $B'_{AD}$ is $\Delta(012) \setminus \{b_{012}\}$. In each of these three cases it is easy to see that $B'_{\tau}$ retracts onto $|I_{\tau}|$.

The case $\tau=\{D\}$ is again special. Here $S_D = \emptyset$ and $B'_D = B_D$ is a four-dimensional subset of $\Delta^4$ missing its barycenter $b_{[4]}$. Pushing radially away from $b_{[4]}$ retracts $B_D$ onto $B_D \cap \partial \Delta^4$, which is easier to describe: we have $$B_D \cap \partial \Delta^4 = \partial \Delta^4 \setminus L,
$$  
where $L$ consists of five line segments joining barycenters of faces in $G_D$ in the cyclic order $
b_{012},b_{024},b_{034},b_{134},b_{123},b_{012}
$. Thus each tetrahedron in $
\partial \Delta^4$ is missing a segment joining the centers of two of its facets. By pushing away from this segment to the remaining two facets and the edge forming the intersection of the two pierced facets, we retract $B_D \cap \partial \Delta^4$ onto the M\"obius band $|I_{D}|$.
\end{document}